\documentclass[10pt]{amsart}
\usepackage{graphicx}
\usepackage{mathrsfs}
\usepackage{epsfig}
\usepackage{amsmath}
\usepackage{amsfonts}
\usepackage{amssymb}
\usepackage{amssymb,amscd}
\usepackage{tikz}
\usepackage{verbatim}
\usepackage{enumerate}
\usepackage{booktabs}
\usepackage[all]{xy}
\usepackage{subfigure}

\newtheorem{thm}{Theorem}[section]

\newtheorem{lem}[thm]{Lemma}
\newtheorem{prop}[thm]{Proposition}
\newtheorem{conj}[thm]{Conjecture}
\theoremstyle{definition}
\newtheorem{defn}[thm]{Definition}

\numberwithin{equation}{section}
\begin{document}

\title[Spherical {CR} uniformization of closed hyperbolic 3-orbifold]{The Menger curve and Spherical {CR} uniformization of a closed hyperbolic 3-orbifold }
\author{Jiming Ma}
\address{School of Mathematical Sciences, Fudan University, Shanghai, 200433, P. R. China}
\email{majiming@fudan.edu.cn}
\author{Baohua Xie}
\address{School of Mathematics, Hunan University, Changsha, 410082, China
}
\email{xiexbh@hnu.edu.cn}



\date{\today}

\thanks{Ma was  supported by NSFC (No.12171092). Xie was supported by NSFC (No.11871202, No.12271148).}
\maketitle

\begin{abstract} 
	Let $G_{6,3}$ be  a hyperbolic polygon-group with boundary the Menger curve. Granier~\cite{Granier}  constructed a  discrete, convex cocompact and  faithful representation $\rho$ of $G_{6,3}$	into $\mathbf{PU}(2,1)$. We show the 3-orbifold at infinity of $\rho(G_{6,3})$  is   a closed hyperbolic 3-orbifold,  with underlying space the 3-sphere  and singular locus  the   $\mathbb{Z}_3$-coned chain-link $C(6,-2)$. 
	This answers the  second part of  Kapovich's Conjecture 10.6 in \cite{Kapovich}, and it also provides the second explicit example of a closed hyperbolic 3-orbifold that admits a uniformizable spherical {CR}-structure after Schwartz's first example in \cite{Schwartz:2003}. 
	
\end{abstract}

\section{Introduction}
% The main difference between them is that the (real) sectional curvature of complex hyperbolic space  is no longer constant, but is
%pinched between $-1$ and $-\frac{1}{4}$, this makes the complex hyperbolic  geometry more difficult to study.

Complex hyperbolic geometry is a  cousin of real hyperbolic geometry, but we currently know very little about it. Let ${\bf H}^2_{\mathbb C}$ be the complex hyperbolic plane. The holomorphic isometry group of ${\bf H}^2_{\mathbb C}$ is $\mathbf{PU}(2,1)$. 
A\emph{ spherical {CR}-structure} on a smooth 3-manifold $M$ is a maximal collection of distinguished charts modeled on the boundary $\partial \mathbf{H}^2_{\mathbb C}$
of  $\mathbf{H}^2_{\mathbb C}$, where coordinate changes are restrictions of transformations from  $\mathbf{PU}(2,1)$.
In other words, a spherical {CR}-structure is a $(G,X)$-structure with $G=\mathbf{PU}(2,1)$ and $X=\mathbb{S}^3$. In contrast to the results on other geometric structures carried on 3-manifolds, there are relatively few examples known of spherical {CR}-structures. A spherical {CR}-structure on a 3-orbifold $M$  is \emph{uniformizable} if it is
obtained as $M=\Gamma\backslash \Omega_{\Gamma}$, where  $\Omega_{\Gamma}\subset \partial \mathbf{H}^2_{\mathbb C}$ is the set of discontinuity of a discrete subgroup  $\Gamma$ acting on $\partial \mathbf{H}^2_{\mathbb C}=\mathbb{S}^3$. The \emph{limit set} $\Lambda_{\Gamma}$ of $\Gamma$ is by definition $\mathbb{S}^3- \Omega_{\Gamma}$.

Goldman and Parker in \cite{GoPa} initiated the study of the deformations of the ideal triangle group in  $\mathbf{PU}(2,1)$. 
They gave an interval contained in the parameter space of complex hyperbolic ideal triangle groups, for points in this interval the corresponding representations are discrete and faithful.
They conjectured that a complex hyperbolic ideal triangle group $\Gamma=\langle I_1, I_2, I_3 \rangle$ is discrete and faithful if and only if $I_1 I_2 I_3$ is not elliptic. Schwartz proved the Goldman-Parker conjecture in \cite{Schwartz:2001ann, Schwartz:2006}. Furthermore, Schwartz analyzed the complex hyperbolic ideal triangle group $\Gamma$ when $I_1 I_2 I_3$ is parabolic and showed  the 3-manifold at infinity of the quotient space ${\bf H}^2_{\mathbb C}/{\Gamma}$ is commensurable with the
Whitehead link complement in the 3-sphere. In other words, the  Whitehead link complement admits uniformizable  spherical {CR}-structure. Schwartz has conjectured the necessary and sufficient condition for a general complex hyperbolic $(p,q,r)$ triangle group $\Delta_{p,q,r}=\langle I_1,I_2,I_3\rangle < \mathbf{PU}(2,1)$ to be a discrete and faithful  representation of an abstract triangle group $T(p,q,r)$ \cite{Schwartz-icm}. Schwartz's conjecture has been proved in a few cases \cite{Deraux:2015, DerauxF:2015, jwx, ParkerWX:2016, ParkerWill:2016}.
These complex hyperbolic triangle groups give more interesting  examples that some cusped hyperbolic 3-manifolds admit  uniformizable spherical {CR}-structures \cite{Acosta:2019, Deraux:2015, DerauxF:2015, jwx, mx, ParkerWill:2016, Schwartz:2001acta}.

These complex hyperbolic triangle groups above  are abstractly commensurable to a surface group or a free group. So these abstract groups are Gromov hyperbolic, with the boundary a circle or a Cantor set. For more details on hyperbolic groups and their boundaries, the reader may refer to  \cite{Gromov}. 
But there are some more complicated groups that can act on ${\bf H}^2_{\mathbb C}$ geometrically.  One example is in \cite{Schwartz:2003}, where  Schwartz
considered an unfaithful  representation of a triangle group into $\mathbf{PU}(2,1)$, the image group  is called $\Delta_{4,4,4;7}$, which is an arithmetic, geometrically finite, discrete subgroup of $\mathbf{PU}(2,1)$. Schwartz determined the 3-manifold at infinity of $\Delta_{4,4,4;7}$ via  a  sophisticated method. It is conjectured that the limit set of Schwartz's group  $\Delta_{4,4,4;7}$ is a Menger curve \cite{Schwartz:2003}.

Let $$G_{6,3}=\langle  a_0, \cdots,  a_5|  a_{i}^{3}=id,  a_{i} a_{i+1}= a_{i+1} a_{i}, i  \in \mathbb{Z}/6\mathbb{Z}\rangle$$	 be  a hyperbolic group with boundary the Menger curve 	\cite{Bourdon}. Recall that the \emph{Menger curve} $\mathcal{K}$ is a one-dimensional locally connected metrizable continuum 
 	without locally separating points, which contain the topological image of any curve. The \emph{standard  Menger curve   in $\mathbb{R}^3$}   can be obtained  as follows: first we subdivide the standard cube $C_{0}=[0,1]^3$ into $3^3$ congruent subcubes; let $C_1$ be the union of these subcubes that intersect the one-skeleton of $C_0$; then we  repeat this process on each subcube again and again to define $C_{n}$; the standard Menger curve in  $\mathbb{R}^3$ is defined  to be the intersection $$\mathcal{M}=\cap^{\infty}_{n=0} C_{n}.$$

Granier  \cite{Granier}  constructed a  discrete, convex-cocompact and  faithful representation $\rho$ of $G_{6,3}$	into $\mathbf{PU}(2,1)$, so the limit set $\Lambda$ of $\rho(G_{6,3})$ is homeomorphic to the boundary of $G_{6,3}$, that is, $\Lambda= \mathcal{K}$ topologically. See Section \ref{sec:rep} for more details on  Granier's representation.

Kapovich made the following conjecture on Granier's representation, see Conjecture 10.6 of \cite{Kapovich}:

\begin{conj}\label{conj:kapovich}
	The Menger curve limit set above is “unknotted” in $\mathbb{S}^3$, i.e.,  the limit set $\Lambda$ of $\rho(G_{6,3})$ is ambient-isotopic to the standard Menger curve $\mathcal{M} \subset \mathbb{R}^3 \subset \mathbb{S}^3$. Furthermore, the quotient 3-dimensional manifold $\Omega/ \rho(G_{6,3})$ is hyperbolic, where  $\Omega$ is the set of discontinuity of $\rho(G_{6,3})$.
\end{conj}

In this paper, we study the topology and geometry of  the 3-orbifold $\Omega/\rho(G_{6,3})$ at infinity of $\rho(G_{6,3})$. We  answer the second  part of Conjecture \ref{conj:kapovich}:
\begin{thm} \label{thm:main}
	The 3-orbifold  $\Omega/\rho(G_{6,3})$  at infinity of $\rho(G_{6,3})$  is   a closed hyperbolic 3-orbifold  $\mathcal{O}$,  with underlying space the 3-sphere  and  singularity locus the   $\mathbb{Z}_3$-coned chain-link $C(6,-2)$.
	
\end{thm}

To the authors' knowledge, the  3-orbifold  $\mathcal{O}$ in Theorem \ref{thm:main} is the second explicit example of a closed hyperbolic  3-orbifold that admits a uniformizable  spherical {CR}-structure after the first example by   Schwartz \cite{Schwartz:2003} nearly twenty years ago.

We prove Theorem \ref{thm:main}  by studying the quotient of the  ideal boundary of  the Dirichlet domain  under the action of the group $\rho(G_{6,3})$.  Granier  described a Dirichlet domain $D$ of $\rho(G_{6,3})$  centered at the fixed point of an  elliptic element of order 6 in $\rho(G_{6,3})$ in her thesis \cite{Granier}.  We will continue to study the topology of $\partial_{\infty}D \cap \Omega$	in $\mathbb{S}^3$. We show that
\begin{thm} \label{thm:unknotdirichlet}
	$\partial_{\infty}D \cap \Omega$	is  a solid torus in the 3-sphere $\partial {\bf H}^2_{\mathbb C}$.
\end{thm}

From Figure \ref{fig:solid-torus} in Section \ref{sec:combinatoricdirichlet}, it seems that 	$\partial_{\infty}D \cap \Omega$	is an unknotted solid torus in the 3-sphere $\partial {\bf H}^2_{\mathbb C}$. That is, the complement of $\partial_{\infty}D \cap \Omega$ in $\partial {\bf H}^2_{\mathbb C}$	is also a solid torus (we do not prove this rigorously), which seems to be strong evidence of the first part of   Conjecture \ref{conj:kapovich}.
%\begin{thm} \label{thm:unknotdirichlet}
%	$D \cap \Omega$	is an unknotted solid torus in the 3-sphere $\partial {\bf H}^2_{\mathbb C}$.
%\end{thm}

\textbf{Outline of the paper}: In Section \ref{sec:background},  we give a well-known background
material on complex hyperbolic geometry. Section \ref{sec:rep}  contains   the matrix representation of the group and the Dirichlet domain constructed by  Granier,  and we will study carefully the combinatorial structure of the ideal boundary of the Dirichlet domain $D$ in Sections \ref{sec:dirichlet} and \ref{sec:combinatoricdirichlet}. In particular, we will prove Theorem \ref{thm:unknotdirichlet} in Section \ref{sec:combinatoricdirichlet}.
We will prove Theorem \ref{thm:main} in Section \ref{sec-3-orbifold} based on results in Sections \ref{sec:dirichlet} and \ref{sec:combinatoricdirichlet}.

\subsection*{Acknowledgements\label{ackowledgements}} We would like to thank the anonymous referees for their helpful comments and constructive suggestions that improved the manuscript. Moreover, one of the referees proposed an  easy way to show the 3-orbifold in Theorem \ref{thm:main} is hyperbolic.
 B. Xie is grateful to  the LMNS (Laboratory of  Mathematics for Nonlinear Science) of Fudan University for its hospitality  during his visit.

\section{Background}\label{sec:background}
In this section, we introduce some background about complex hyperbolic geometry. Almost all facts
stated here can be found in the book of Goldman \cite{Go} and in   \cite{dpp:2016}.

Let $\mathbb{C}^{2,1}$ be  the $3$-dimensional complex vector space consisting of $3$-tuples
\begin{equation*}
Z=\left(
\begin{array}{c}
z_1 \\
z_2 \\
z_3 \\
\end{array}
\right)\in \mathbb{C}^{3}
\end{equation*} endowed with a Hermitian form $\langle \cdot, \cdot\rangle$, which has signature $(2,1)$.

A vector $Z$ is said to be negative (respectively null, positive) if and only if the Hermitian form $\langle Z, Z\rangle$
is negative (respectively null, positive).
Let $\mathbb{P}:\mathbb{C}^{2,1}\backslash\{0\}\rightarrow \mathbb{CP}^2$ be the standard projection map.
\emph{Complex hyperbolic plane ${\bf H}^2_{\mathbb C}$} is defined to be the subset of $\mathbb{P}(\mathbb{C}^{2,1} \backslash\{0\})$
consisting of negative lines in $\mathbb{C}^{2,1}$.
The boundary of ${\bf H}^2_{\mathbb C}$ is the subset $\partial{\bf H}^2_{\mathbb C}$ of $\mathbb{P}(\mathbb{C}^{2,1}\backslash\{0\})$
consisting of null lines in $\mathbb{C}^{2,1}$.

In this paper, we will use two different models of complex hyperbolic plane. There are two different Hermitian matrices $J$
which give different Hermitian forms on $\mathbb{C}^{2,1}$. Let $Z,W$ be the column vectors $(z_1,z_2,z_3)^{t}$ and
$(w_1,w_2,w_3)^{t}$ respectively. The first Hermitian form is defined to be
\begin{equation*}
\langle Z,W\rangle=-z_1\overline{w_1}+z_2\overline{w_2}+z_3\overline{w_3}.
\end{equation*}
It is given by the Hermitian matrix $J_1$:
\begin{equation*}
J_1=\left[
\begin{array}{ccc}
-1 & 0 & 0 \\
0 & 1 & 0 \\
0 & 0 & 1 \\
\end{array}
\right].
\end{equation*}
Note that this Hermitian form agrees with the one given in \cite{Granier}.

The second Hermitian form is defined to be
\begin{equation*}
\langle Z,W\rangle=z_1\overline{w_3}+z_2\overline{w_2}+z_3\overline{w_1}.
\end{equation*}
It is given by the Hermitian matrix $J_2$:
\begin{equation*}
J_2=\left[
\begin{array}{ccc}
0 & 0 & 1\\
0 & 1 & 0 \\
1 & 0 & 0 \\
\end{array}
\right].
\end{equation*}

We define the first model of complex hyperbolic plane by taking $z_1=1$ in column vector
$Z=(z_1,z_2,z_3)^{t}\in \mathbb{C}^{2,1}$ for the first Hermitian form.
We then have

\begin{equation*}
{\bf H}^2_{\mathbb C}=\left\{\left(
\begin{array}{c}
1\\
z_1 \\
z_2 \\
\end{array}
\right)\in \mathbb{CP}^{2}\bigg| |z_1|^{2}+|z_2|^{2}<1\right\}.
\end{equation*}
This forms the unit ball model of complex hyperbolic plane. The boundary $\partial{\bf H}^2_{\mathbb C}$ is the sphere $\mathbb{S}^3$ given by
$$|z_1|^{2}+|z_2|^{2}=1.$$
Almost all calculations in Sections \ref{sec:dirichlet} and \ref{sec:combinatoricdirichlet} will be done in the ball model. For the convenience to drawing pictures, we also consider the Siegel model.
We obtain the Siegel model of complex hyperbolic plane by taking $z_3=1$ in column vector $Z=(z_1,z_2,z_3)^{t}\in \mathbb{C}^{2,1}$
for the second Hermitian form. That is,

\begin{equation*}
{\bf H}^2_{\mathbb C}=\left\{\left(
\begin{array}{c}
z_1 \\
z_2 \\
1
\end{array}
\right)\in \mathbb{CP}^{2}\bigg| 2 {\rm Re}(z_1)+|z_2|^{2}<0\right\}.
\end{equation*}
Its boundary $\partial{\bf H}^2_{\mathbb C}$ is
\begin{equation*}
\partial{\bf H}^2_{\mathbb C}=\left\{\left(-(|z|^2+is)/2,z,1)^t|(z,s)\in \mathbb{C}\times \mathbb{R}\right)\}\cup \{\infty=(1,0,0)^t\right\}.
\end{equation*}
It is also the one point compactification of the 3-dimensional Heisenberg group $\mathbb{C}\times \mathbb{R}$, with group law
\begin{equation*}
[z,t]\cdot[w,s]=[z+w,t+s+2{\rm Im}(z\overline{w})].
\end{equation*}

To move from the ball model to the Siegel model we introduce the Cayley transformation
\begin{equation*}
	C=\left[
	\begin{array}{ccc}
		\frac{1}{\sqrt{2}} & 0 & \frac{1}{\sqrt{2}}\\
		0 & 1 & 0 \\
		-\frac{1}{\sqrt{2}} & 0 & \frac{1}{\sqrt{2}}\\
	\end{array}
	\right].
\end{equation*}

The group of linear isometries preserving the Hermitian form J is a non-compact group isomorphic to $\mathbf{ U}(2,1)$ (with respect to J).
We denote   $\mathbf{ PU}(2,1)$ by the group of holomorphic isometries of ${\bf H}^2_{\mathbb C}$,  which is the projectivization of the
unitary group $\mathbf{ U}(2,1)$. We will often consider matrices in the group $\mathbf{ SU}(2,1)$ instead of elements of $\mathbf{ PU}(2,1)$.
Every element of $\mathbf{ PU}(2,1)$  admits exactly three lifts to the group  $\mathbf{ SU}(2,1)$ of unitary matrices  for J of determinant one.

Up to scaling, ${\bf H}^2_{\mathbb C}$ carries a unique $\mathbf{ U}(2,1)$-invariant Riemannian metric with curvature between $-1$ and $-1/4$.
The metric information we will use is the following distance formula:
\begin{equation}\label{eq:distance}
\cosh\left(\frac{d(u,v)}{2}\right)=\frac{|\langle \mathbf{ u},\mathbf{ v}\rangle|}{\sqrt{\langle \mathbf{ u},\mathbf{ u}\rangle \langle \mathbf{ v},\mathbf{ v}\rangle}},
\end{equation} where $\mathbf{ u}, \mathbf{ v}$ denote lifts of $u, v$ to $\mathbb{C}^3$.

If $u,v\in \mathbb{C}^3$, we define the \emph{Hermitian cross product} of $u$ and  $v$, denoted by $u\boxtimes v$, as the Euclidean cross product of the
vectors $u^*J$ and $v^*J$, where $J$ is the matrix defining the  Hermitian form,    $u^*$ and $ v^*$ are the conjugate transpose vectors of $u$ and $v$.
If $u$ and $v$ are collinear, then $u\boxtimes v=0$.  If not, then $u\boxtimes v$ spans their Hermitian orthogonal complement and
$\langle u\boxtimes v, u\rangle=0, \langle u\boxtimes v, v\rangle=0$.

\subsection{Totally geodesic subspace}

Given a positive vector $v \in \mathbb{C}^{2,1}$, its orthogonal complement $v^{\perp}=\{u\in \mathbb{C}^3: \langle  v, u\rangle=0\}$ is a two
dimensional subspace on which the Hermitian form restricts to a form with signature $(1,1)$. The set of negative lines in $v^{\perp}$ is
then a copy of ${\bf H}^1_{\mathbb C}$, naturally isometric to the Poincar\'e disk.

\begin{defn} The submanifold  of ${\bf H}^2_{\mathbb C}$ given by $v^{\perp}$ is called a \emph{complex geodesic}. The vector
	$v$ is called a \emph{polar} to the complex geodesic $v^{\perp}$.
\end{defn}

Given a vector $v$ with $\langle v,v\rangle=1$, we consider the isometry of  $\mathbb{C}^3$ given by
\begin{equation*}
R_{v,\zeta}(x)=x+(\zeta-1)\langle x,v\rangle v
\end{equation*} where $\zeta$ is a complex number of absolute value one.  It is easy to see that $R_{v,\zeta}(x)$ preserves the Hermitian
inner product,  fixes  the vectors in $v^{\perp}$, and rotates the normal direction by an angle $\theta$, where $\zeta=e^{i\theta}$.

\begin{defn}The isometry $R_{v,\zeta}$ is called a \emph{complex reflection} with mirror $v^{\perp}$.
\end{defn}

Another type of totally geodesic subspace is given in the standard ball model as the set of points with real coordinates, which is the fixed point set of the isometry $(x_1,x_2)\rightarrow(\overline{x_1},\overline{x_2})$. It is simply a copy of the real hyperbolic plane ${\bf H}^2_{\mathbb R}$. This submanifold of ${\bf H}^2_{\mathbb C}$  is often called a \emph{Langrangian plane}.

The boundary at infinity of a complex geodesic is called \emph{$\mathbb{C}$-circle}, and the boundary at infinity of a copy of ${\bf H}^2_{\mathbb R}$ is called  an
\emph{$\mathbb{R}$-circle}. The group $\mathbf{ PU}(2,1)$ acts transitively on each kind of subspaces.

\subsection{Bisectors and their intersections}
Note that there are no totally geodesic real hypersurfaces in complex hyperbolic space. The Dirichlet polyhedra are bounded by bisectors, which are hypersurfaces
equidistant from two given points. Their geometric structure and complicated intersection patterns have been analyzed in great detail in \cite{Go}. We will review some of the
results which will be needed in the paper.

\begin{defn}
	The \emph{bisector} between two distinct points  $p_0$ and $p_1$ in  ${\bf H}^2_{\mathbb C}$ is the set of points that are equidistant from  $p_0$ and $p_1$:
	\begin{equation}\label{eq:bisector1}
	\mathcal{B}(p_0,p_1)=\{u\in {\bf H}^2_{\mathbb C}: d(u,p_0)=d(u, p_1)\}.
	\end{equation}
\end{defn}
We  denote  $\mathbf{p_0}$ and $\mathbf{p_1}$  lifts  of $p_0$ and $p_1$ to $\mathbb C^3$.
In view of equation (\ref{eq:distance}), if we normalize two vectors $\mathbf{p_0}$ and $\mathbf{p_1}$  so that $\langle \mathbf{p_0}, \mathbf{p_0}\rangle=\langle \mathbf{p_1},\mathbf{p_1} \rangle$, the equation (\ref{eq:bisector1})
of the bisector then becomes simply
\begin{equation}\label{eq:bisector2}
|\langle u, \mathbf{p_0}\rangle|=|\langle u,\mathbf{p_1} \rangle|.
\end{equation}

A bisector  in  ${\bf H}^2_{\mathbb C}$ is a smooth codimension one real hypersurface diffeomorphic  to a 3-ball, but it is
not totally geodesic.  The \emph{spinal sphere} of the bisector  $\mathcal{B}(p_0, p_1)$ is the ideal boundary of it on $\partial {\bf H}^2_{\mathbb C}$.
 The \emph{complex spine} of the bisector  $\mathcal{B}(p_0, p_1)$ is by definition the complex geodesic that contains
$p_0$ and $p_1$. The \emph{real spine} is a real geodesic in the complex spine that is equidistant between $p_0$ and $p_1$.
There is a natural extension of the real spine to projective space, given by the (not necessarily negative) vectors  in $Span_\mathbb{C}(\mathbf{p_0},\mathbf{p_1})$ satisfying
equation (\ref{eq:bisector2}).  We
call this the \emph{extended real spine}  of the bisector.

Although bisector $\mathcal{B}$ is not totally geodesic, it can be described in two different ways in terms of a foliation by totally geodesic subspaces. $\mathcal{B}$ is the preimage of the real spine
under  the orthogonal projection onto the complex spine; each complex line that is the  preimage of a point of the real spine is  called a \emph{complex slice} of $\mathcal{B}$. $\mathcal{B}$ is also the union of all Lagrangian planes that contain the real spine; each Lagrangian plane is called a \emph{real slice} of $\mathcal{B}$.

%There are two kinds of maximal totally geodesic submanifolds contained in a given bisector $\mathcal{B}$. The complex ones are called  \emph{complex slices} of $\mathcal{B}$; they are
%the  preimages of points of the real spine under the orthogonal projection onto the complex spine.  The totally real submanifolds that are contained in $\mathcal{B}$ are precisely
%those containing the real spine, and are called \emph{real slices} of $\mathcal{B}$.

We will  describe bisectors by giving two points on their extended real spines.  In fact, we can describe a real geodesic
as the projectivization  of a totally real $2$-dimensional  subspace of $\mathbb{C}^{2,1}$, i.e. we take two vectors $u$ and $v$ in $\mathbb{C}^{3}$ with
$\langle u, v\rangle\in \mathbb{R}$, and consider their real span.  The simplest way to guarantee that the span really yields a geodesic in  ${\bf H}^2_{\mathbb C}$
is to require moreover that $v$ and $u$ form a Lorentz basis, i.e. $\langle v, v\rangle=-1, \langle u, u\rangle=1$ and $\langle v, u\rangle=0$.

Let $\Sigma_1$ and $\Sigma_2$ be the complex spines of the bisectors $\mathcal{B}_1$ and $  \mathcal{B}_2$, and let $\sigma_1$ and $\sigma_2$ be their real spines. Then $\Sigma_1$ and $\Sigma_2$ coincide, or they intersect at a single point. Note that if $\Sigma_1$ and $\Sigma_2$ intersect outside the real spines $\sigma_1$ and $\sigma_2$, then $\mathcal{B}_1\cap  \mathcal{B}_2$ can be written as the equidistant locus from three points in general position.

First, we assume that the bisectors $\mathcal{B}_1$ and $\mathcal{B}_2$ have the same complex spine.
\begin{defn}\label{def:cospinal}
	The bisectors $\mathcal{B}_1$ and $\mathcal{B}_2$ are called \emph{cospinal} if and only if their complex spines $\Sigma_1$ and  $\Sigma_2$ coincide.
\end{defn}
In this case, it follows from the slice decomposition that $\mathcal{B}_1\cap\mathcal{B}_2$ is non-empty if and only if their real spines $\sigma_1$ and  $\sigma_2$ intersect at a point $p \in {\bf H}^2_{\mathbb C}$. Moreover, $\mathcal{B}_1\cap\mathcal{B}_2$ consists of a complex geodesic $\mathcal{S}$, namely  the complex geodesic orthogonal to $\Sigma_1=\Sigma_2$ through the point $p$.

Furthermore, we can describe the intersection of $\mathcal{S}$  with another bisector $\mathcal{B}_3(p_0,p_3)$ as follows, see \cite{Deraux:2006}.

We can choose a basis $\{v_1, v_2\}$ for $\mathcal{S}$, with $\langle v_1,v_2\rangle=0$, $\langle v_1,v_1\rangle=-1$, and $\langle v_2,v_2\rangle=1$.
Then the vectors in $\mathcal{S}$ can be parameterized as
\begin{equation}\label{eq:complex-geodesic}
\mathcal{S}=\{v_1+zv_2: |z|<1\}.
\end{equation}

The intersection $\mathcal{S}\cap \mathcal{B}_3(p_0,p_3) $ with a third bisector has  an equation of the form
\begin{equation}\label{eq:complexline-bisector}
|\langle v_1+zv_2,p_0\rangle|=|\langle v_1+zv_2,p_3\rangle|,
\end{equation}
which is a circle (or Euclidean line) in the $z$-plane.

In particular, the intersection of $\mathcal{S}$ with a number of half spaces in ${\bf H}^2_{\mathbb C}$ is bounded by circles or lines in the unit disk. One should be aware
that such an intersection need not be connected in general, since the arcs of circles bounding it are not necessarily geodesic.

Now suppose $\Sigma_1$ and $\Sigma_2$ are distinct complex spines of $\mathcal{B}_1$ and $\mathcal{B}_2$ respectively. We call $\mathcal{B}_1$ and $\mathcal{B}_2$ \emph{coequidistant} if and only if $\Sigma_1$ and $\Sigma_2$ intersect outside the real
spines. The following result
is  due to Giraud, see \cite{Go}, which is crucial to the study of Dirichlet domain.

\begin{prop}\label{prop:giraud}Let $p_0$, $p_1$ and $p_2$ be distinct points in ${\bf H}^2_{\mathbb C}$, not all contained in a complex line. When it is non-empty, the intersection
	$\mathcal{B}(p_0,p_1)\cap\mathcal{B}(p_0,p_2)$ is a (non-totally geodesic) smooth disk. Moreover, it is contained in precisely three bisectors, namely
	$\mathcal{B}(p_0,p_1), \mathcal{B}(p_0,p_2)$ and $\mathcal{B}(p_1,p_2)$.
\end{prop}

\begin{defn}The intersection of two coequidistant bisectors with distinct complex spines is called \emph{Giraud disk}.
\end{defn}

In order to find the combinatorics of the Dirichlet domain,  we must determine  the intersection of the bisectors.  Unlike in real hyperbolic space, the intersection of two bisectors is not necessarily connected, see  \cite{Go} for the details. Due to Proposition \ref{prop:giraud}, we choose the bisectors bounding a Dirichlet domain $D$ to form a very special family of bisectors; namely, they are all equidistant from a given point $p_0$. 

We will review a convenient way to  parametrize a Giraud disk,  see  \cite{DerauxF:2015,  dpp:2016}. Consider two  coequidistant bisectors $\mathcal{B}_1(p_0,p_1)$ and $\mathcal{B}_2(p_0,p_2)$,
which we assume not to be cospinal.

Let  $\mathbf{p}_j$ denote a lift of $p_j$ to $\mathbb{C}^3$.  We may assume that  the three square norms
$\langle\mathbf{p}_j,\mathbf{p}_j \rangle$ are equal up to the rescaling. We also assume  $\langle\mathbf{p}_0,\mathbf{p}_1 \rangle$ and $\langle\mathbf{p}_0,\mathbf{p}_2 \rangle$ are
real and positive.

Define $\widetilde{v_j}=\mathbf{p}_0-\mathbf{p}_j$ and $\widetilde{w_j}=i(\mathbf{p}_0+\mathbf{p}_j)$ for $j=1,2$.  We can normalize these to unit vectors
$v_j=\widetilde{v_j}/\sqrt{-\langle\widetilde{v_j},\widetilde{v_j} \rangle}$ and $w_j=\widetilde{w_j}/\sqrt{\langle\widetilde{w_j},\widetilde{w_j} \rangle}$. Then $\widetilde{v_j}$ is the midpoint of the geodesic segment between $p_0$ and $p_j$. The extended real spine of $\mathcal{B}_j(p_0,p_j)$ 
can be parametrized as the (not necessarily negative) vectors of the form $w_j+t v_j$  ($t\in \mathbb{R}$).
 Thus the intersection $\mathcal{B}_1\cap\mathcal{B}_2$
is given by negative vectors of the form
\begin{equation*}\label{eq:Gdisk}
V(t_1,t_2)=(w_1+t_1v_1)\boxtimes(w_2+t_2v_2),
\end{equation*}
with $t_1,t_2\in \mathbb{R}$.  Its extension to projective space will be called \emph{Giraud torus}, i.e., the vector $V(t_1,t_2)$ is no longer required to have a negative square norm, we often denote it by  $\widehat{\mathcal{B}_1} \cap\widehat{\mathcal{B}_2}$.
With the parameterization of the extended real spine of $\mathcal{B}_j$, it is evident that point $v_j$ is missing. But the projectivization of the orthogonal complement of $v_j$($j=1,2$) does not intersect ${\bf H}^2_{\mathbb C}$.
  Given three points $p_0$, $p_1$ and $p_2$, it is easy to determine whether the intersection $\mathcal{B}_1(p_0,p_1)\cap\mathcal{B}_2(p_0,p_2)$ is empty or not by finding a sample point. Note that the condition $\langle V(t_1,t_2),V(t_1,t_2)\rangle$ is negative is equivalent to
\begin{equation}\label{eq:Gdisk-negative}
\det\left[
\begin{array}{cc}
\langle w_1+t_1v_1,w_1+t_1v_1\rangle&\langle w_1+t_1v_1,w_2+t_2v_2\rangle\\
\langle w_2+t_2v_2,w_1+t_1v_1\rangle&\langle w_2+t_2v_2,w_2+t_2v_2\rangle
\end{array}
\right]>0.
\end{equation}

\section{The representation and the group}\label{sec:rep}

In her Ph.D thesis \cite{Granier},  Granier constructed a convex-compact representation  $\rho$ of the polygon-group $G_{6,3}$ in $\mathbf{ PU}(2,1)$. She constructed a Dirichlet domain
for $\rho(G_{6,3})$ and proved that $\rho(G_{6,3})$ is discrete by using Poincar\'{e} polyhedron theorem in ${\bf H}^2_{\mathbb C}$.
\begin{defn}Fix two natural numbers $p\geq 5$ and $q\geq 3$. The \emph{polygon-group} is defined as follows
	\begin{equation}
G_{p,q}=	\langle a_{0},a_{1},\ldots a_{p-1}|a_{i}^{q}=[a_i,a_{i+1}]=id, i   \in \mathbb{Z}/p\mathbb{Z}\rangle.
	\end{equation}
\end{defn}

Let 
\begin{equation}
H_{p,q}=	\langle a_{0},r|a_{0}^{q}=r^{p}=[a_0,ra_{0}r^{-1}]=id \rangle
\end{equation}
be another group.  Let $\phi: H_{p,q}\longrightarrow \mathbb{Z}_p$ be the group homomorphism defined by $\phi(a_0)=0$ and $\phi(r)=1$. Then $G_{p,q}$ is the kernel of $\phi$. This shows that $G_{p,q}$ is a normal subgroup of $H_{p,q}$ of index $p$.  We owe this observation to one of the referees.

We now review the representation $\rho$ of $G_{6,3}$  in \cite{Granier}. In fact, Granier gave the  representation $\rho:H_{6,3} \rightarrow \mathbf{PU}(2,1)$.
We will use the first model for the complex hyperbolic plane in Section \ref{sec:background}. Consider a regular right-angled $p$-gon $\mathcal{P}$ in a well-chosen real hyperbolic plane
${\bf H}^2_{\mathbb R}\subset {\bf H}^2_{\mathbb C}$ with vertices $x_j$ for $j=0,1\cdots p-1$. The center of $\mathcal{P}$ is located at $o=(1,0,0)^t$. The lift of $x_j$ is given as
\begin{equation*}
x_j=\left(
\begin{array}{c}
1 \\
s \cdot \cos(\frac{2 j\pi}{p})  \\
s \cdot  \sin(\frac{2 j\pi}{p})
\end{array}
\right),
\end{equation*} where $$s=\frac{\sqrt{2\cos^2(\frac{2 \pi}{p})+2\cos(\frac{2 \pi}{p})}}{1+\cos(\frac{2 \pi}{p})},$$ which is the Euclidean distance between the origin $o$ and $x_j$. The Euclidean distance $s$ is related to the complex hyperbolic distance $d(o,x_j)$ by
\begin{equation*}
s=\tanh\left(\frac{d(o,x_j)}{2}\right), 
\end{equation*} from which it follows
\begin{equation*}
d(o,x_j)=2 \operatorname{arccosh} \left(\dfrac{1+\cos(2\pi/p)}{\sin(2\pi/p)}\right).
\end{equation*}

For $0\leq j\leq p-1$, the geodesic side $[x_j,x_{j+1}]$ of $\mathcal{P}$ determines a complex geodesic $C_j$ with polar vector $e_j=x_j\boxtimes x_{j+1}$. For example, one can get
\begin{equation*}
e_0=\left(
\begin{array}{c}
\sqrt{2\cos^{2}(2\pi/p)+2\cos(2\pi/p)} \\
1+\cos(2\pi/p) \\
\sin(2\pi/p)
\end{array}
\right).
\end{equation*}

Let $\gamma_j$ be the complex reflection of order $q$ with mirror $C_j$ and let
\begin{equation*}
R_p=\left[
\begin{array}{ccc}
1 & 0 & 0\\
0 & \cos(2\pi/p) & -\sin(2\pi/p) \\
0 & \sin(2\pi/p) & \cos(2\pi/p)\\
\end{array}
\right].
\end{equation*}

We see that $R_p(x_j)=x_{j+1}$ and the polygon $\mathcal{P}$ is preserved by the rotation $R_p$. It follows that $R_p(C_j)=C_{j+1}$ and $\gamma_{j}=R_p^{j}\gamma_0 R_p^{-j}$.
We define $\rho$ by the map  $\rho(a_j)=\gamma_j$  and $\rho(r)=R_{p}$.

When $p=6$ and $q=3$, we have

\begin{equation*}
R_6=\left[
\begin{array}{ccc}
1 & 0 & 0\\
0 & \frac{1}{2} & -\frac{\sqrt{3}}{2} \\
0 & \frac{\sqrt{3}}{2} &\frac{1}{2}\\
\end{array}
\right]
\end{equation*}

and
\begin{equation*}
\gamma_0=\left[
\begin{array}{ccc}
\frac{5}{2}-i\frac{\sqrt{3}}{2} & -\frac{3 \sqrt{6}}{4}+i\frac{3\sqrt{2}}{4} & -\frac{3 \sqrt{2}}{4}+i\frac{\sqrt{6}}{4}\\
\frac{3 \sqrt{6}}{4}-i\frac{3\sqrt{2}}{4} & -\frac{5}{4}+i\frac{3\sqrt{3}}{4} & -\frac{3 \sqrt{3}}{4}+\frac{3}{4}i \\
\frac{3 \sqrt{2}}{4}-i\frac{\sqrt{6}}{4} & -\frac{3 \sqrt{3}}{4}+\frac{3}{4}i & \frac{1}{4}+i\frac{\sqrt{3}}{4}\\
\end{array}
\right].
\end{equation*}

Since the polar vectors $e_j$ and $e_{j+1}$ are orthogonal, we have $\gamma_i\gamma_{i+1}=\gamma_{i+1}\gamma_i$.   Subsequently, we often write $\Gamma$ instead of $\rho(G_{6,3})$.
%An easy computation shows that $\gamma_i\gamma_{i+1}=\gamma_{i+1}\gamma_i$ for $0 \leq i \leq 5$.   Subsequently, we often write $\Gamma$ instead of $\rho(G_{6,3})$.

\section{Dirichlet domain of the group}\label{sec:dirichlet}

Given any group $G$ acting on  ${\bf H}^2_{\mathbb C}$, the Dirichlet domain centered at $p_0 \in {\bf H}^2_{\mathbb C}$ is by definition

\begin{equation*}
D=D_{G}=\{u\in {\bf H}^2_{\mathbb C}: d(u,p_0)\leqslant d(u, \gamma(p_0)), \forall \gamma\in G\}.
\end{equation*}

The group $G$ acts discretely if and only if $D$ has  nonempty interior, and in that case,  $D_{G}$ is a fundamental domain for $G$ as long as
no element of the group fixes the point $p_0$. If $p_0$ is fixed by some non-trivial element of $G$, then one only gets a fundamental domain modulo the action of the finite group that fixes $p_0$.

Note that $R_6$ is a regular elliptic element of order 6, with isolated fixed point at $o=(1,0,0)^{t}$. The point $o$ will be chosen as the center of the Dirichlet domain of $\Gamma$.

One wishes the Dirichlet domain of $\Gamma$ has only finitely many faces that are on bisectors, so that we can first consider the \emph{partial Dirichlet domain}
\begin{equation*}
D_{S}=\{u\in {\bf H}^2_{\mathbb C}: d(u,p_0)\leqslant d(u, \gamma(p_0)), \forall \gamma\in S\}
\end{equation*} for some finite subset $S\subset \Gamma$. Then $D_{S}$ will be  a priori domain larger than $D_{\Gamma}$ by taking into account only the faces coming from $S$ rather than all of $\Gamma$.

In order to ensure that $D_{S}$ has side pairings, the generating set $S$ should be symmetric, that is, $S$ is closed under the operation of taking inverses in the group. For the group
$\Gamma=\rho(G_{6,3})$, one can guess a reasonable candidates for the set $S$, which is the following set of 24 group elements:
\begin{equation*}
S=\{\gamma_i,\gamma_i^{-1},\gamma_i\gamma_{i+1}^{-1}, \gamma_{i+1}\gamma_i^{-1}, i\in \mathbb{Z}/6\mathbb{Z} \}.
\end{equation*}

Consider the partial Dirichlet domain $D_{S}$. For a given $\gamma\in S$, note that the intersection $D_{S}\cap \mathcal{B}(o,\gamma(o))$ may be very complicated. When $D_{S}\cap \mathcal{B}(o,\gamma(o))$ has a nonempty interior in $\mathcal{B}(o,\gamma(o))$, the face of $D_{S}$ associated to the element $\gamma$ is given by the connected components of $D_{S}\cap \mathcal{B}(o,\gamma(o))$, which is often denoted by $b(o,\gamma(o))$.

The main result of \cite{Granier} is the following.

\begin{thm}The Dirichlet domain $D_{\Gamma}$ centered at $o$ is equal to  $D_{S}$. In particular, $D_{\Gamma}$ has precisely 24 faces, namely the faces of $D_{\Gamma}$ associated to the elements of $S$.
\end{thm}

In order to prove that $D_{S}$ is a fundamental domain $D_{\Gamma}$ of $\Gamma$, one should start by determining the  precise combinatorics of $D_{S}$, then check the conditions of the  Poincar\'e polyhedron theorem.

%Let $\Gamma$ be an index 6 subgroups of $\Gamma_{6,3}$ is generated by $\gamma_0,\ldots,\gamma_5$
%\begin{equation*}
%\Gamma=\{\gamma_i |[\gamma_i,\gamma_{i+1}]=id, i\in \mathbb{Z}/6\mathbb{Z} \}.
%\end{equation*}

Note that $D_{S}$ is a fundamental domain of $\rho(G_{6,3})$, but it is not a fundamental domain of  $\rho(H_{6,3})$. In fact, since the  element $ R_6$  fixes the center of $D_{S}$  and preserve the $D_{S}$ from the construction, $D_{S}$ is a fundamental domain for
the coset decomposition of  $\rho(H_{6,3})$ into left cosets of the cyclic group $\langle R_6\rangle$.

For convenience, we write $\mathcal{B}_k$,  $i=0, 1,\ldots,11$,  for the bisectors bounding $D_{S}$,  and $b_k$, $i=0,1,\ldots,11$,  for the intersection
$\mathcal{B}_k\cap D_{S}$.  See Table 1 for the notations. We also write $\mathcal{B}_{\overline{k}},  b_{\overline{k}}$ for the bisectors and faces that associate to the inverse of the elements  which associate to the $\mathcal{B}_k$ and $b_k$. For example, $\mathcal{B}_0=\mathcal{B}(o,\gamma_0(o))$ and $\mathcal{B}_{\overline{0}}=\mathcal{B}(o,\gamma_0^{-1}(o))$.
For each $k$, we denote  $\overline{\mathcal{B}_k}, \overline{\mathcal{B}_{\overline{k}}}$ and $\overline{b_k}, \overline{b_{\overline{k}}}$ be their closures in $\overline{{\bf H}^2_{\mathbb C}}={\bf H}^2_{\mathbb C}\cup\partial{\bf H}^2_{\mathbb C}$.

We next describe the symmetry of $D_{S}$. From the construction, $D_{S}$ is $R_6$-invariant. It has at most 2 isometry types of faces, that is, each face is isometric to the face on
$\mathcal{B}_0$ or  $\mathcal{B}_6$.  In \cite{Granier}, Granier also observed that there are two complex reflections  preserving $D_{S}$.

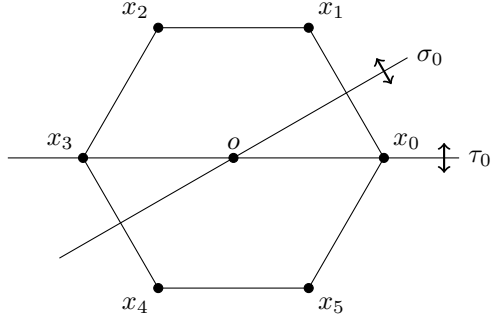
\begin{figure}
	\begin{center}
		\begin{tikzpicture}[scale=2]
		
		\coordinate (O) at (0,0);
		\coordinate (x0) at (1,0);
		\coordinate (x1) at ([shift={+(120:1cm)}] x0);
		\coordinate (x2) at ([shift={+(-1,0)}] x1);
		\coordinate (x3) at ([shift={+(240:1cm)}] x2);
		\coordinate (x4) at ([shift={+(300:1cm)}] x3);
		\coordinate (x5) at ([xshift=1cm, yshift=0] x4);
		
		\draw (x0)node [above right]{$x_0$}
		--(x1)node [above right]{$x_1$}
		--(x2)node [above left]{$x_2$}
		--(x3)node [above left]{$x_3$}
		--(x4)node [below left]{$x_4$}
		--(x5)node [below right]{$x_5$}
		--cycle;
		
		\draw (-1.5,0) -- (1.5, 0);
		\draw [thick,<->] (1.4,0.1) -- (1.4, -0.1);\node [right] at (1.5,0){$\tau_0$};
		\draw (-1.157,-2/3) -- (1.157, 2/3);
		\draw [thick,<->] (0.95, 0.663953) -- (1.05, 0.490748);\node [right] at (1.157, 2/3){$\sigma_0$};
		\node [above] at (0,0){$o$};
		\draw [fill] (x0) circle [radius=.03];
		\draw [fill] (x1) circle [radius=.03];
		\draw [fill] (x2) circle [radius=.03];
		\draw [fill] (x3) circle [radius=.03];
		\draw [fill] (x4) circle [radius=.03];
		\draw [fill] (x5) circle [radius=.03];
		\draw [fill] (O) circle [radius=.03];

		\end{tikzpicture}
		
	\end{center}
	\caption{The polygon  $\mathcal{P}$ and its symmetry.}
	\label{figure:symmetry-P}
\end{figure}

Let $\tau_0, \sigma_0 : {\mathbb{C}}^3 \longrightarrow {\mathbb{C}}^3$ be given as follows:
\begin{equation*}
\tau_0 : \left(
\begin{array}{c}
z_1 \\
z_2 \\
z_3 \\
\end{array}
\right) \longmapsto
\left(
\begin{array}{c}
z_1 \\
z_2 \\
-z_3 \\
\end{array}
\right)
\end{equation*}
and
\begin{equation*}
\sigma_0 : \left(
\begin{array}{c}
z_1 \\
z_2 \\
z_3 \\
\end{array}
\right) \longmapsto
\left(
\begin{array}{c}
z_1 \\
z_2\cos(\frac{2\pi}{6})+z_3\sin(\frac{2\pi}{6}) \\
z_2\sin(\frac{2\pi}{6})-z_3\cos(\frac{2\pi}{6})  \\
\end{array}
\right).
\end{equation*} Then we have
\begin{itemize}
	
	\item $\tau_0 $ fixes the vertices $x_0$ and $x_3$, which also  interchanges two pair of vertices $\{x_1,x_5\}$ and $\{x_2,x_4\}$;
	
	\item $\tau_0 \gamma_{k} \tau_0=\gamma_{5-k}$ for $k=0\ldots,5$, therefore,  $\tau_0(\mathcal{B}_k)=\mathcal{B}_{5-k}$ for $k=0\ldots,5$;   $\tau_0(\mathcal{B}_k)=\mathcal{B}_{\overline{16-k}}$ for  $k=6\ldots,10$  and  $\tau_0(\mathcal{B}_{11})=\mathcal{B}_{\overline{11}}$;
	
	\item $\sigma_0 $ interchanges three pair of vertices $\{x_0,x_1\}, \{x_2,x_5\}$ and $\{x_3,x_4\}$;
	
	\item $\sigma_0 \gamma_{0} \sigma_0=\gamma_{0}$,  $\sigma_0 \gamma_{k} \sigma_0=\gamma_{6-k}$ for $k=1,2,3,4,5$, therefore, $\sigma_0 (\mathcal{B}_{0})=\mathcal{B}_{0}$;  $\sigma_0 (\mathcal{B}_{k}) =\mathcal{B}_{6-k}$ for $k=1,2,3,4,5$; $\sigma_0 (\mathcal{B}_{\overline{0}})=\mathcal{B}_{\overline{0}}$,  $\sigma_0 (\mathcal{B}_{\overline{k}}) =\mathcal{B}_{\overline{6-k}}$ for $k=1,2,3,4,5$; $\sigma_0 (\mathcal{B}_{k}) =\mathcal{B}_{\overline{17-k}}$ for $k=6,\ldots,11$. See  Figure \ref{figure:symmetry-P}.

\end{itemize}

Define $\tau_j=R_6^{j}\tau_0R_6^{-j}$ and $\sigma_j=R_6^{j}\sigma_0R_6^{-j}$. It is easy to see that these isometries also preserve $D_{S}$.

\begin{table}[!htbp]
	\centering
	\caption{The notations for the bisectors  associated with $D_{S}$.}
	\label{table1}
	
	\begin{tabular}{ccc}
		\toprule
		\textbf{Element of $S$ } & \textbf{Bisector  } \\
		\midrule
		$\gamma_0$&$\mathcal{B}_0=\mathcal{B}(o,\gamma_0(o))$\\
		$\gamma_1$&$\mathcal{B}_1=R_6(\mathcal{B}_0)$\\
		$\gamma_2$&$\mathcal{B}_2=R_6^2(\mathcal{B}_0)$\\
		$\gamma_3$&$\mathcal{B}_3=R_6^3(\mathcal{B}_0)$\\
		$\gamma_4$&$\mathcal{B}_4=R_6^4(\mathcal{B}_0)$\\
		$\gamma_5$&$\mathcal{B}_5=R_6^5(\mathcal{B}_0)$\\
		$\gamma_0\gamma_1^{-1}$&$\mathcal{B}_6=\mathcal{B}(o,\gamma_0\gamma_1^{-1}(o))$\\
		$\gamma_1\gamma_2^{-1}$&$\mathcal{B}_7=R_6(\mathcal{B}_6)$\\
		$\gamma_2\gamma_3^{-1}$&$\mathcal{B}_8=R_6^2(\mathcal{B}_6)$\\
		$\gamma_3\gamma_4^{-1}$&$\mathcal{B}_9=R_6^3(\mathcal{B}_6)$\\
		$\gamma_4\gamma_5^{-1}$&$\mathcal{B}_{10}=R_6^4(\mathcal{B}_6)$\\
		$\gamma_5\gamma_0^{-1}$&$\mathcal{B}_{11}=R_6^5(\mathcal{B}_6)$\\
		%$\gamma_{0}^{-1}$&$\mathcal{B}_{\overline{0}}=\mathcal{B}(o,\gamma_0^{-1}(o))$&$b_{\overline{0}}$\\
		%$\gamma_{1}^{-1}$&$\mathcal{B}_{\overline{1}}=R_6(\mathcal{B}_{\overline{0}})$&$b_{\overline{1}}$\\
		%$\gamma_{2}^{-1}$&$\mathcal{B}_{\overline{2}}=R_6^2(\mathcal{B}_{\overline{0}})$&$b_{\overline{2}}$\\
		%$\gamma_{3}^{-1}$&$\mathcal{B}_{\overline{3}}=R_6^3(\mathcal{B}_{\overline{0}})$&$b_{\overline{3}}$\\
		%$\gamma_{4}^{-1}$&$\mathcal{B}_{\overline{4}}=R_6^4(\mathcal{B}_{\overline{0}})$&$b_{\overline{4}}$\\
		%$\gamma_{5}^{-1}$&$\mathcal{B}_{\overline{5}}=R_6^5(\mathcal{B}_{\overline{0}})$&$b_{\overline{5}}$\\
		%$\gamma_1\gamma_0^{-1}$&$\mathcal{B}_{\overline{6}}=\mathcal{B}(o,\gamma_1\gamma_0^{-1}(o))$&$b_{\overline{6}}$\\
		%$\gamma_2\gamma_1^{-1}$&$\mathcal{B}_{\overline{7}}=R_6(\mathcal{B}_{\overline{6}})$&$b_{\overline{7}}$\\
		%$\gamma_3\gamma_2^{-1}$&$\mathcal{B}_{\overline{8}}=R_6^2(\mathcal{B}_{\overline{6}})$&$b_{\overline{8}}$\\
		%$\gamma_4\gamma_3^{-1}$&$\mathcal{B}_{\overline{9}}=R_6^3(\mathcal{B}_{\overline{6}})$&$b_{\overline{9}}$\\
		%$\gamma_5\gamma_4^{-1}$&$\mathcal{B}_{\overline{10}}=R_6^4(\mathcal{B}_{\overline{6}})$&$b_{\overline{10}}$\\
		%$\gamma_0\gamma_5^{-1}$&$\mathcal{B}_{\overline{11}}=R_6^5(\mathcal{B}_{\overline{6}})$&$b_{\overline{11}}$\\
		\bottomrule
		
	\end{tabular}
	
\end{table}

\subsection{The combinatorics of $D_{S}$}

We now describe the  combinatorics of $D_{S}$ following \cite{Granier} in details. By the symmetry of  $D_{S}$, it is enough to determine the combinatorics of two faces of
$D_{S}$, namely,  $\mathcal{B}_0\cap D_{S}$ and  $\mathcal{B}_6\cap D_{S}$. Since any two of these 24 bisectors bounding $D_{S}$ are coequidistant, their pairwise intersection is diffeomorphic to a disk, which is either a Giraud disk or a complex geodesic.

The following two lemmas give the details of the intersections of $\mathcal{B}_0$ with the other 23 bisectors.

\begin{lem}$\mathcal{B}_0$ intersects exactly  9 bisectors of the  23 other bisectors, that is, $\mathcal{B}_{\overline{0}}$,  $\mathcal{B}_1$, $\mathcal{B}_{\overline{1}}$, $\mathcal{B}_5$, $\mathcal{B}_{\overline{5}}$, $\mathcal{B}_6$, $\mathcal{B}_{\overline{6}}$, $\mathcal{B}_{11}$ and $\mathcal{B}_{\overline{11}}$. The intersection $\mathcal{B}_0\cap \mathcal{B}_{\overline{0}}$ is a complex geodesic and the other corresponding intersections are all Giraud disks.
\end{lem}

\begin{lem}
	$\mathcal{B}_0\cap \mathcal{B}_{\overline{1}} \cap D_{S}$, $\mathcal{B}_0\cap \mathcal{B}_{\overline{5}} \cap D_{S}$, $\mathcal{B}_0\cap \mathcal{B}_{\overline{6}} \cap D_{S}$ and $\mathcal{B}_0\cap \mathcal{B}_{11} \cap D_{S}$ are empty set.
\end{lem}

A $2$-face of $D_{S}$ is the $2$-dimensional intersection of two faces of $D_{S}$.    The precise combinatorics of each $2$-face of $b_0$ has been studied in detail in  \cite{Granier}.

\begin{figure}[htbp]
	\centering
	\begin{minipage}[t]{0.48\textwidth}
		\centering
		\includegraphics[width=5cm]{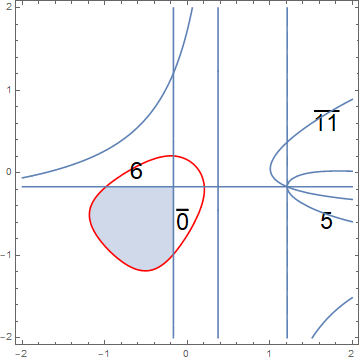}
		\caption{The $2$-face of $\overline{b_0}$ on the Giraud disk $\overline{\mathcal{B}_0}\cap \overline{\mathcal{B}_{1}}$, which is bounded by the red Jordan curve. Here we write $k$ for $\mathcal{B}_{k} $ and $\overline{k}$ for $\mathcal{B}_{\overline{k}}$.}\label{fig:2-face-B0-and-B1}
	\end{minipage}
	\begin{minipage}[t]{0.48\textwidth}
		\centering
		\includegraphics[width=5cm]{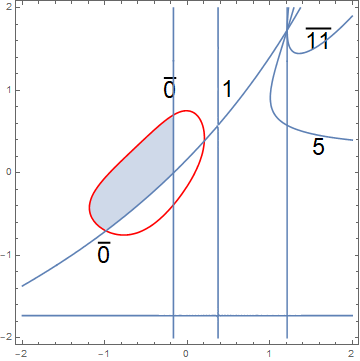}
		\caption{The $2$-face of $\overline{b_0}$ on the Giraud disk $\overline{\mathcal{B}_0}\cap \overline{\mathcal{B}_{6}}$.}
		\label{fig:2-face-B0-and-B6}
	\end{minipage}
\end{figure}

\begin{figure}[htbp]
	\centering
	\begin{minipage}[t]{0.48\textwidth}
		\centering
		\includegraphics[width=5cm]{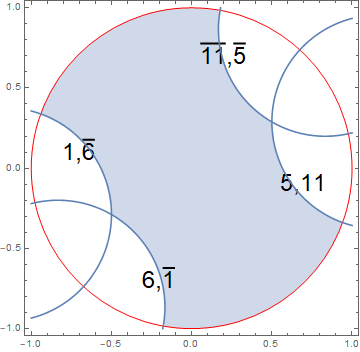}
		\caption{The $2$-face of $\overline{b_0}$ on the  complex geodesic $\overline{\mathcal{B}_0}\cap \overline{\mathcal{B}_{\overline{0}}}. $ }\label{complex-geodesic }
		\label{fig:2-face-B0-and-barB1}
	\end{minipage}
	\begin{minipage}[t]{0.48\textwidth}
		\centering
		\includegraphics[width=5cm]{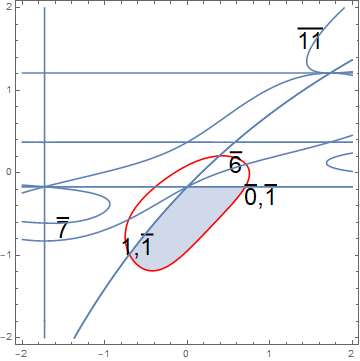}
		\caption{The $2$-face of $\overline{b_6}$ on the Giraud disk $\overline{\mathcal{B}_6}\cap \overline{\mathcal{B}_{0}} $.}
		\label{fig:2-face-B6-and-barB0}
	\end{minipage}
\end{figure}

\begin{prop}\label{prop:b0} The closure $\overline{b_0}$ of $b_0$ in ${\bf H}^2_{\mathbb C}\cup\partial {\bf H}^2_{\mathbb C}$ has precisely six $2$-faces, five finite ones and one on the spinal sphere associated to  $\mathcal{B}_0$.
	\begin{itemize}
		\item  The finite $2$-faces on the (closure of the) Giraud disks $\overline{\mathcal{B}_0}\cap \overline{\mathcal{B}_{1}} $ and $\overline{\mathcal{B}_0}\cap \overline{\mathcal{B}_{5}} $ are  topological triangles, see Figure 	\ref{fig:2-face-B0-and-B1}. In particular, the second $2$-face is the image of the first $2$-face under the action of the isometry $\sigma_0$.
		\item  The finite $2$-faces on the (closure of the) Giraud disks $\overline{\mathcal{B}_0}\cap \overline{\mathcal{B}_{6}} $ and $\overline{\mathcal{B}_0}\cap \overline{\mathcal{B}_{\overline{11}}} $ are  topological triangles, see Figure 	\ref{fig:2-face-B0-and-B6}. In particular, the second $2$-face is the image of the first $2$-face under the action of the isometry $\sigma_0$.
		\item  The finite $2$-face on the (closure of the) complex geodesic $\overline{\mathcal{B}_0}\cap \overline{\mathcal{B}_{\overline{0}}} $ is a topological hexagon, see Figure 	\ref{fig:2-face-B0-and-barB1}.
		
		\item  The  $2$-face on the spinal sphere $\partial\mathcal{B}_0$ is a topological  hexagon, see Figure 	\ref{fig:schematicb0} and Figure \ref{figure:sidepair}.
		
	\end{itemize}
	
\end{prop}

\begin{lem}$\mathcal{B}_6$ intersects exactly  7 bisectors of the other 23 bisectors,  that is, $\mathcal{B}_{0}$,   $\mathcal{B}_{\overline{0}}$, $\mathcal{B}_1$, $\mathcal{B}_{\overline{1}}$, $\mathcal{B}_{\overline{6}}$, $\mathcal{B}_{\overline{7}}$ and $\mathcal{B}_{\overline{11}}$.  The corresponding intersections are all Giraud disks.
\end{lem}

\begin{lem}
	$\mathcal{B}_6\cap \mathcal{B}_{\overline{0}} \cap D_{S}$, $\mathcal{B}_6\cap \mathcal{B}_{1} \cap D_{S}$, $\mathcal{B}_6\cap \mathcal{B}_{\overline{6}} \cap D_{S}$,  $\mathcal{B}_6\cap \mathcal{B}_{\overline{7}} \cap D_{S}$ and  $\mathcal{B}_6\cap \mathcal{B}_{\overline{11}} \cap D_{S}$ are empty.
\end{lem}

The precise combinatorics of each $2$-face of $b_6$ has also been studied in detail in  \cite{Granier}.

\begin{prop}\label{prop:b6} The closure $\overline{b_6}$ of $b_6$ in ${\bf H}^2_{\mathbb C}\cup\partial {\bf H}^2_{\mathbb C}$ has precisely three $2$-faces, two finite ones and one on the spinal sphere associated to  $\mathcal{B}_6$.
	\begin{itemize}
		\item  The finite $2$-faces on the (closure of the) Giraud disks $\overline{\mathcal{B}_6}\cap \overline{\mathcal{B}_{0}} $ and $\overline{\mathcal{B}_6}\cap \overline{\mathcal{B}_{\overline{1}}} $ are  topological triangles, see Figure \ref{fig:2-face-B6-and-barB0}. In particular, the second $2$-face is the image of the first $2$-face under the action of the isometry $\sigma_0$.

		\item  The  $2$-face on the spinal sphere $\partial\mathcal{B}_6$ is a   bigon, see Figure 	\ref{fig:schematicbarb6} and Figure \ref{figure:sidepair}.
		
	\end{itemize}
	
\end{prop}

By applying a version of Poincar\'{e} polyhedron theorem in the complex hyperbolic plane  as stated for example in \cite{ParkerWill:2016}, \cite{dpp:2016} or \cite{Mostow:1980}, the main result obtained in \cite{Granier} can be stated as follows.
\begin{thm}
	Let $D_{S}$ be defined  as above, then $D_{S}$ is a fundamental domain for  $\rho(G_{6,3})$.
	Moreover, $\Gamma=\rho(G_{6,3})$ is discrete and has the presentation
	$$
	\langle \gamma_0,\ldots,\gamma_5| \gamma^3_j=id, \gamma_j\gamma_{j+1}=\gamma_{j+1}\gamma_j, j\in \mathbb{Z}/6\mathbb{Z}\rangle.
	$$
\end{thm}

\begin{figure}[htbp]
	\centering
	\begin{minipage}[t]{0.49\textwidth}
		\centering
		\begin{tikzpicture}
		\draw [red, thick] (1,0.5) circle [radius=2.5];
		\draw [fill] (0.5,0.5) circle [radius=0.1];
		\draw [fill] (1.5,0.5) circle [radius=0.1];
		\draw [fill] (3.5,0.5) circle [radius=0.1];
		\draw [fill] (-1.5,0.5) circle [radius=0.1];
		\draw [thick] (0.5,0.5) -- (-1.5,0.5);
		\draw [thick] (1.5,0.5) -- (3.5,0.5);
		\draw [thick] (0,2.79) -- (0.5,0.5);
		\draw [fill] (0,2.79) circle [radius=0.1];
		\draw [thick] (0,-1.79) -- (0.5,0.5);
		\draw [fill] (0,-1.79) circle [radius=0.1];		
		\draw [thick] (2,2.79) -- (1.5,0.5);
		\draw [fill] (2,2.79) circle [radius=0.1];
		\draw [thick] (2,-1.79) -- (1.5,0.5);
		\draw [fill] (2,-1.79) circle [radius=0.1];
		
		\node [left] at (-0.5,1.2) {$\mathbf{1}$};
		\node [left] at (-0.5,-0.2) {$\mathbf{6}$};
		\node [left] at (-0.5,0.5) {$\mathbf{\overline{1}}$};		
		\node [left] at (3,1.2) {$\mathbf{\overline{11}}$};
		\node [left] at (3,-0.2) {$\mathbf{5}$};
		\node [left] at (3,0.5) {$\mathbf{\overline{5}}$};
		\node [above] at (1,0.5) {$\mathbf{\overline{0}}$};
		\node [above] at (0.4,1.5) {$\mathbf{\overline{6}}$};
		\node [below] at (0.4,-0.5) {$\mathbf{\overline{1}}$};
		\node [above] at (1.6,1.5) {$\mathbf{\overline{5}}$};
		\node [below] at (1.6,-0.5) {$\mathbf{11}$};
		\end{tikzpicture}
		\caption{A schematic view of the combinatorics of the face of $\overline{b_0}$.}
		\label{fig:schematicb0}
	\end{minipage}
	\begin{minipage}[t]{0.49\textwidth}
		\centering
		\begin{tikzpicture}
		\draw [red, thick] (1,0.5) circle [radius=2.5];
		\draw [thick] (1,0.5) -- (1,3);
		\draw [thick] (1,0.5) -- (1,-2);
		\draw [fill] (1,0.5) circle [radius=0.1];
		\draw [fill] (1,3) circle [radius=0.1];
		\draw [fill] (1,-2) circle [radius=0.1];
		\node [left] at (0.1,0.5) {$\mathbf{0}$};
		\node [right] at (1.9,0.5) {$\mathbf{\overline{1}}$};
		\node [right] at (1,1.7) {$\mathbf{\overline{0}}$};
		\node [right] at (1,-0.7) {$\mathbf{1}$};
		\end{tikzpicture}
		\caption{A schematic view of the combinatorics of the face of $\overline{b_6}$.}
		\label{fig:schematicbarb6}
	\end{minipage}
\end{figure}

\section{The combinatorics at infinity of  the Dirichlet domain}\label{sec:combinatoricdirichlet}

The ideal boundary of a set $X\subset \mathbf{H}^2_{\mathbb C} $ will be
denoted by $\partial_{\infty}X\subset \partial\mathbf{H}^2_{\mathbb C} $.
Let $\Omega$ be   the set of discontinuity of the discrete subgroup  $\rho(G_{6,3})$ acting on $\partial \mathbf{H}^2_{\mathbb C}=\mathbb{S}^3$. In this section, we show	$\partial_{\infty}D_{S} \cap \Omega$	is a solid torus in the 3-sphere $\partial {\bf H}^2_{\mathbb C}$. 
The idea is to consider the intersection with
$\partial {\bf H}^2_{\mathbb C}$ of the fundamental domain $D_{S}$ for the action on $ {\bf H}^2_{\mathbb C}$.  The main result in this section is Proposition  \ref{prop:octagon}, which is also the key for the proof of Theorem  \ref{thm:main}.

\begin{figure}[htbp]
		\centering
		\includegraphics[width=14cm]{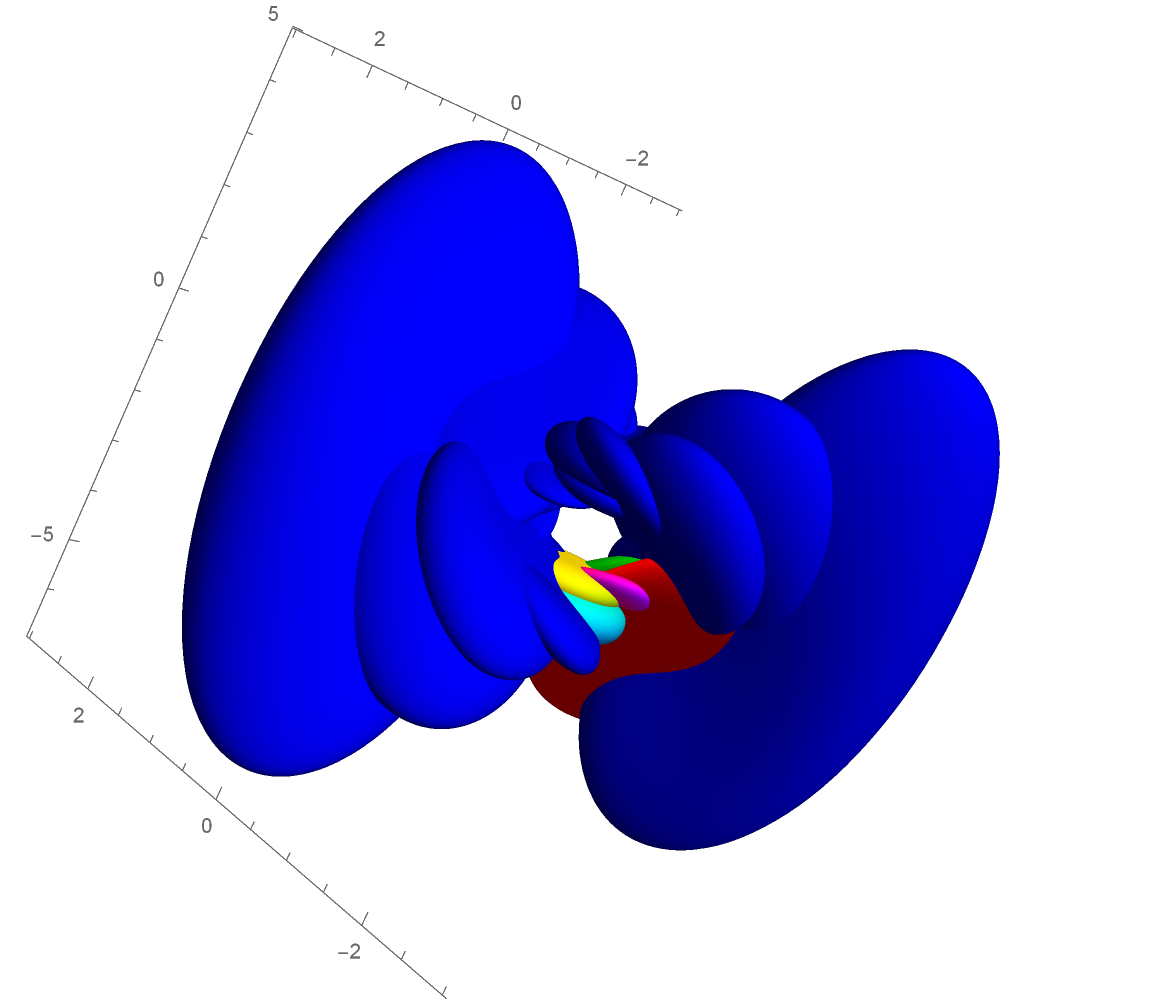}
		\caption{The solid torus $T$ is drawn on the boundary of Siegel model. Where the red (resp. green, yellow, cyan, purple) sphere is the spinal sphere of $\gamma_{0}$ (resp. $\gamma^{-1}_{0}$, $\gamma^{-1}_{1}$, $\gamma_{1}$, $\gamma_{6}$). }\label{fig:solid-torus}	
\end{figure}

In what follows, the set $\partial_{\infty}D_{S}$ will be denoted by $T$. Note that $T$ is bounded by 24 pieces of spinal spheres. A realistic view of $T$ is given in Figure \ref{fig:solid-torus}, which is drawn in the boundary of Siegel model. So  $T$ is contained in the unbounded region bounded by the torus  in Figure \ref{fig:solid-torus}. From the analysis of the combinatorics of $D_{S}$ given in the previous section, the combinatorial structure of $\partial T$ can be seen in a schematic picture given in Figure 	\ref{figure:sidepair}. The picture is obtained by putting together the incidence  information for each $2$-face on these spinal spheres. One should keep this picture in mind for the gluing of these 24 faces.

\begin{prop}  $\partial T$ is a torus.
\end{prop}
\begin{proof}
	From the precise combinatorics of the 2-faces of $b_0$ and $b_6$
given in Propositions \ref{prop:b0}, \ref{prop:b6}, we can get  Figure \ref{figure:sidepair} by gluing together  these 24 faces.
The curves colored red and green are glued together, respectively.
Thus we see that $\partial T$ is a torus.

%In fact, one can see that $T$ is a solid torus by Proposition \ref{prop:octagon}.
% This is  clear from  Figure \ref{figure:sidepair}.   
\end{proof}

 The torus $\partial T$ divides $\mathbb{S}^3$ into two parts: outside (containing $\infty$) and inside. The inside may be a solid torus and the outside a knot complement, or vice-versa, or both sides may be solid tori. J. W. Alexander's theorem tells us that $\partial T$ in $\mathbb{S}^3$ bounds a solid torus on at least one side, see \cite{Rolfsen}. In fact, in our case it seems  $\partial T$ bounds a solid torus on both side  in $\mathbb{S}^3$. 
We just show $T$ is a solid torus by producing an explicit simple closed curve
which  bounds a disk in $T$.

%We will show this by produce two explicit simple closed curves with intersection number one on $\partial T$, onebounds a disk in one side of $\mathbb{S}^3-\partial T$.

The main goal in our analysis is to produce an explicit embedded disk in $T$ whose boundary is the red curve on the up and down side of Figure  \ref{figure:sidepair}. This disk will cut $T$
into a 3-ball, then we can get the fundamental group of $\Omega/\Gamma$ from the gluing maps.

\begin{figure}[htbp]
	\centering
	\includegraphics[width=14cm]{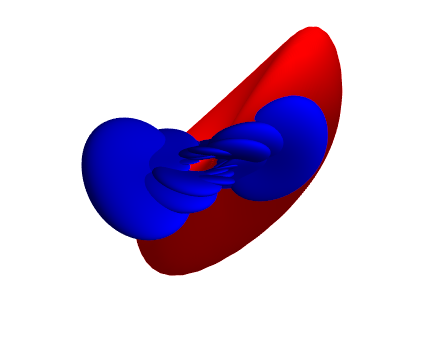}
	\caption{The solid torus $T$ and the cutting disk on the spinal sphere of  the bisector $\mathcal{C}$.  We draw all the 24 spinal spheres for the set $S$ in blue, and the spinal sphere of $\mathcal{C}$ in Proposition  \ref{prop:cutbisector} in red. $T$ intersects this red sphere in two disks, one of which is our  topological octagon  $\mathcal{E}$ in Proposition \ref{prop:octagon}, which in turn implies $T$ is a solid torus.}
	\label{fig::cuttingdisk}	
\end{figure}

%\begin{figure}[htbp]
%	\centering
%	\begin{minipage}[t]{0.48\textwidth}
%		\centering
%		\includegraphics[width=5cm]{torus.png}
%		\caption{The solid torus $T$ is drawn on the boundary of Siegel model. }\label{fig:solid-torus}
%	\end{minipage}
%	\begin{minipage}[t]{0.48\textwidth}
%		\centering
%		\includegraphics[width=5cm]{torus-disk.png}
%		\caption{The solid torus $T$ and the cut disk on the spinal sphere of  the bisector $\mathcal{C}$.}
%		\label{fig:cuttingdisk}
%	\end{minipage}
%\end{figure}

Let $\mathcal{C}$ be the bisector $\mathcal{B}(\gamma_1\gamma_0^{-1}(o),\gamma_0\gamma_1^{-1}(o))$. Then one can show that

%\begin{prop} \label{prop:cutbisector}The bisector $\mathcal{C}$ intersects with the eights %bisectors $\mathcal{B}_{\overline{0}}$,  $\mathcal{B}_{11}$, $\mathcal{B}_{5}$, %$\mathcal{B}_{4}$, $\mathcal{B}_{3}$, $\mathcal{B}_{2}$, $\mathcal{B}_{\overline{7}}$, %$\mathcal{B}_{\overline{1}}$.
	
%\end{prop}

%\begin{proof}
%\end{proof}
\begin{prop}\label{prop:cutbisector} The bisector $\mathcal{C}$ intersects with the 24 bisectors bounding $D_{S}$.	
\end{prop}

\begin{proof} 
	
We show that the intersection of the bisector $\mathcal{C}$ with each bisector of the 24 bisectors is a smooth disk. The argument used here can be found in Appendix of \cite{DerauxF:2015}.
	In order to show that $\mathcal{C}\cap \mathcal{B}_{\overline{0}}$ is a disk, we consider the parametrization of the Giraud torus $\widehat{\mathcal{C}}\cap \widehat{\mathcal{B}_{\overline{0}}}$. We only need to exhibit a single point $p_{\overline{0}}\in {\bf H}^2_{\mathbb C}$ inside the Giraud torus. For example, the point
	$p_{\overline{0}}$ has the following form
	\begin{eqnarray*}
		p_{\overline{0}}&=&(\gamma_0\gamma_1^{-1}(o)-\gamma_1\gamma_0^{-1}(o))\boxtimes(o-\gamma_0^{-1}(o))\\
		&=&\left(\frac{-3\sqrt{3}-9i}{2},\frac{-3-3\sqrt{3}i}{2\sqrt{2}},\frac{-3\sqrt{3}-9i}{2\sqrt{2}}\right)
	\end{eqnarray*}
	does the job, since $\langle p_{\overline{0}}, p_{\overline{0}}\rangle=-9$.
	
	For the remainder cases, we just list the points and the corresponding  Giraud tori. See Table \ref{tab:tvertices}.
\end{proof}
	
%	\begin{table}
%		\begin{tabular}{c|c}
%			Giraud torus & the point                   \\\hline
%			$\widehat{\mathcal{C}}\cap \widehat{\mathcal{B}_{11}}$ & $p_{11}=\left(\frac{-9\sqrt{3}+9i}{2},\frac{-9+4\sqrt{3}i}{\sqrt{2}},-3\sqrt{\frac{3}{2}}\right)$ \\
%			$\widehat{\mathcal{C}}\cap \widehat{\mathcal{B}_{5}}$  & $p_{5}=\left(\frac{9-9\sqrt{3}+(9\sqrt{3}-9)i}{2},\frac{19\sqrt{3}-39+i(33-15\sqrt{3})}{2\sqrt{2}},\frac{3-3\sqrt{3}+i(3\sqrt{3}-9)}{2\sqrt{2}}\right)$ \\
%			$\widehat{\mathcal{C}}\cap \widehat{\mathcal{B}_{4}}$   & $p_{4}=\left(-\frac{3\sqrt{3}+9i}{2},\frac{3-3\sqrt{3}i}{2\sqrt{2}},\frac{3\sqrt{3}-9i}{2\sqrt{2}}\right)$  \\
%			$\widehat{\mathcal{C}}\cap \widehat{\mathcal{B}_{3}}$   & $p_{3}=\left(-\frac{3\sqrt{3}+9i}{2},\frac{3-3\sqrt{3}i}{2\sqrt{2}},\frac{3\sqrt{3}-9i}{2\sqrt{2}}\right)$  \\
%			$\widehat{\mathcal{C}}\cap \widehat{\mathcal{B}_{2}}$   & $p_{2}=\left(-9+\frac{9i}{2}+\frac{9\sqrt{3}}{2},\frac{9\sqrt{3}-18+i(-3-2\sqrt{3})}{2\sqrt{2}},\frac{-45+24\sqrt{3}+i(12-3\sqrt{3})}{2\sqrt{2}}\right)$ \\
%			$\widehat{\mathcal{C}}\cap \widehat{\mathcal{B}_{\overline{7}}}$   & $p_{\overline{7}} =\left(\frac{3\sqrt{3}i-27+12\sqrt{3}}{2},\frac{4\sqrt{3}-9+i(6-7\sqrt{3})}{2\sqrt{2}},\frac{-54+27\sqrt{3}+3i}{2\sqrt{2}}\right)$ \\
%			$\widehat{\mathcal{C}}\cap \widehat{\mathcal{B}_{\overline{1}}}$&$p_{\overline{1}}=\left(-\frac{3\sqrt{3}-9i}{2},\frac{-3-3\sqrt{3}i}{2\sqrt{2}},\frac{-3\sqrt{3}-9i}{2\sqrt{2}}\right)$\\
%		\end{tabular}
%		\normalsize
%		\caption{The point chosen inside each Giraud torus.}\label{tab:tvertices}
%	\end{table}

	\begin{table}
	\begin{tabular}{c|c}
	Point &      	Giraud torus              \\\hline
		$\left(\frac{-3\sqrt{3}+9i}{2},\frac{-3+3\sqrt{3}i}{2\sqrt{2}},\frac{-3\sqrt{3}+9i}{2\sqrt{2}}\right)$ &$\widehat{\mathcal{C}}\cap \widehat{\mathcal{B}_{0}}$, $\widehat{\mathcal{C}}\cap \widehat{\mathcal{B}_{1}}$, $\widehat{\mathcal{C}}\cap \widehat{\mathcal{B}_{6}}$, $\widehat{\mathcal{C}}\cap \widehat{\mathcal{B}_{\overline{0}}}$, $\widehat{\mathcal{C}}\cap\widehat{\mathcal{B}_{\overline{1}}}$, $\widehat{\mathcal{C}}\cap\widehat{\mathcal{B}_{\overline{6}}}$\\
		$\left(\frac{3\sqrt{3}-9i}{2},\frac{-3+3\sqrt{3}i}{2\sqrt{2}},\frac{-3\sqrt{3}+9i}{2\sqrt{2}}\right)$ &$\widehat{\mathcal{C}}\cap \widehat{\mathcal{B}_{3}}$, $\widehat{\mathcal{C}}\cap \widehat{\mathcal{B}_{4}}$, $\widehat{\mathcal{C}}\cap \widehat{\mathcal{B}_{9}}$, $\widehat{\mathcal{C}}\cap \widehat{\mathcal{B}_{\overline{3}}}$, $\widehat{\mathcal{C}}\cap\widehat{\mathcal{B}_{\overline{4}}}$, $\widehat{\mathcal{C}}\cap\widehat{\mathcal{B}_{\overline{9}}}$\\
		$\left(\frac{9\sqrt{3}-81i}{14},\frac{-9+27\sqrt{3}i}{14\sqrt{2}},\frac{9\sqrt{3}-81i}{14\sqrt{2}}\right)$ &$\widehat{\mathcal{C}}\cap \widehat{\mathcal{B}_{2}}$, $\widehat{\mathcal{C}}\cap \widehat{\mathcal{B}_{7}}$,  $\widehat{\mathcal{C}}\cap \widehat{\mathcal{B}_{\overline{2}}}$, $\widehat{\mathcal{C}}\cap\widehat{\mathcal{B}_{\overline{7}}}$\\
		$\left(\frac{18\sqrt{3}-27i}{7},\frac{18\sqrt{2}-9\sqrt{6}i}{7},0\right)$ &$\widehat{\mathcal{C}}\cap \widehat{\mathcal{B}_{5}}$, $\widehat{\mathcal{C}}\cap \widehat{\mathcal{B}_{11}}$,  $\widehat{\mathcal{C}}\cap \widehat{\mathcal{B}_{\overline{5}}}$, $\widehat{\mathcal{C}}\cap\widehat{\mathcal{B}_{\overline{11}}}$\\
		$\left(\frac{-111\sqrt{3}-477i}{62},\frac{111+159\sqrt{3}i}{31\sqrt{2}},0\right)$ &$\widehat{\mathcal{C}}\cap \widehat{\mathcal{B}_{8}}$,  $\widehat{\mathcal{C}}\cap\widehat{\mathcal{B}_{\overline{8}}}$\\
		$\left(\frac{111\sqrt{3}-477i}{62},\frac{111-159\sqrt{3}i}{62\sqrt{2}},\frac{-111\sqrt{3}+477i}{62\sqrt{2}}\right)$ &$\widehat{\mathcal{C}}\cap \widehat{\mathcal{B}_{10}}$,  $\widehat{\mathcal{C}}\cap\widehat{\mathcal{B}_{\overline{10}}}$\\
	\end{tabular}
	\normalsize
	\caption{The point chosen inside each Giraud torus.}\label{tab:tvertices}
\end{table}

\begin{figure}[htbp]
	\centering
	\includegraphics[totalheight=3in]{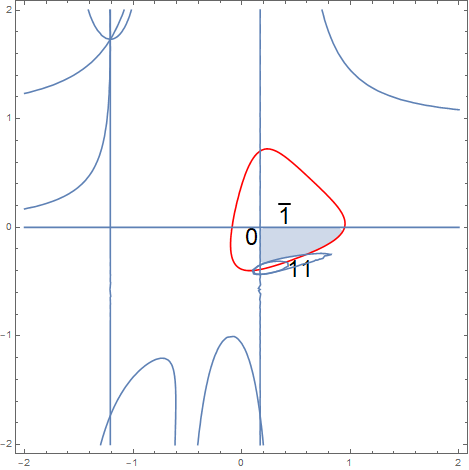}
	\caption{The combinatorics of $\mathcal{C}\cap \mathcal{B}_{\overline{0}}\cap D_{S}$. The ideal boundary of the  intersection is an arc with endpoints on $\partial_{\infty}\mathcal{B}_{\overline{1}}$ and $\partial_{\infty}\mathcal{B}_{11}$. }\label{fig:C-Bm0}	
\end{figure}
By using the techniques of Section \ref{sec:background}, we can study the intersection of $\mathcal{C}$ with each face of $D_{S}$.  For
example, the combinatorics of $\mathcal{C}\cap \mathcal{B}_{0}\cap D_{S}$ is shown in Figure \ref{fig:C-Bm0}. Then one can get the combinatorics of
the intersection of the sphere $\partial_{\infty}\mathcal{C}$ with
the solid torus $T$.  It can be seen that the interior of this intersection has two components, each of which is a topological octagon.
Let $\mathcal{E}$ be one of them.  We will describe this  octagon explicitly. The thick gray curve in  Figure \ref{figure:sidepair} in Section \ref{sec-3-orbifold} is a schematic view of the boundary of  $\mathcal{E}$. From which, it can be seen that $\partial \mathcal{E}$ passes through the spinal spheres of  $\mathcal{B}_{\overline{0}}$,  $\mathcal{B}_{11}$, $\mathcal{B}_{5}$, $\mathcal{B}_{4}$, $\mathcal{B}_{3}$, $\mathcal{B}_{2}$, $\mathcal{B}_{\overline{7}}$, and $\mathcal{B}_{\overline{1}}$ circularly.

We study the intersection of the closure of the bisectors  with $\partial_{\infty}\mathcal{C}$ by parametrizing $\partial_{\infty}\mathcal{C}$.  In order to make the equation of these intersections having simple forms, we will choose the coordinates for $ {\bf H}^2_{\mathbb C}$ such that the midpoint of $[\gamma_1\gamma_0^{-1}(o),\gamma_0\gamma_1^{-1}(o)]$ is at the
origin. The coordinate transformation matrix is given by
\begin{equation*}
P=\left[
\begin{array}{ccc}
4\sqrt{\frac{2}{5}} & 0 & -3\sqrt{\frac{3}{5}}i\\
\frac{3}{2}\sqrt{\frac{3}{5}}&-\frac{\sqrt{3}}{2}i & -2\sqrt{\frac{2}{5}}i \\
\frac{9}{2\sqrt{5}}&\frac{i}{2} & -2\sqrt{\frac{6}{5}}i\\
\end{array}
\right].
\end{equation*}

We work in $\mathbb{C}^2$, with affine coordinates $u_1=\frac{z_1}{z_0}$, $u_2=\frac{z_2}{z_0}$, where the $z_j$ denote coordinates in the new Lorentz basis. The ball coordinates for
$\gamma_1\gamma_0^{-1}(o)$ and $\gamma_0\gamma_1^{-1}(o)$ are given by $(\pm \sqrt{\frac{3}{5}},0)$, and the bisector $\mathcal{C}$  has a very simple equation, namely
$${\rm Re}(u_1)=0.$$  So the bisector $\mathcal{C}$  can simply be though of as the unit ball in $\mathbb{R}^3$,

\begin{equation*}
\mathcal{C}=\{(it_3,t_1+it_2)\in \mathbb{C}^2 |t_i\in \mathbb{R}, t_1^2+ t_2^2+ t_3^2<1\},
\end{equation*}
and $\partial_{\infty}\mathcal{C}$ is the unit sphere.

The equation for the intersection of $\mathcal{C}$ with the bisector $\mathcal{B}(o,g(o))$ for some $g$ has the form
\begin{equation*}
|\langle Z,P^{-1}\left(p_0\right)\rangle|=|\langle Z,P^{-1}\left(g(p_0)\right)\rangle|, 
\end{equation*}
where one takes $Z=(1,it_3,t_1+it_2)$.

\begin{table}[!htbp]
	\caption{The equations for the intersections of the eight bisectors with the bisector $\mathcal{C}$.}
	\centering
	\begin{tabular}{c|c}
		\toprule
		
		$\mathcal{B}_{\overline{0}}$& $7+12t_1^2+12t_2^2+6\sqrt{5}t_3-5t_3^2+9\sqrt{6}t_2+\sqrt{30}t_2t_3-5\sqrt{2}t_1+3\sqrt{10}t_1t_3$\\
		$\mathcal{B}_{11}$&$113+78t_1^2+78t_2^2-60\sqrt{5}t_3+35t_3^2+20\sqrt{2}t_1-12\sqrt{10}t_1t_3+76\sqrt{6}t_2-20\sqrt{30}t_2t_3$\\
		$\mathcal{B}_{5}$&$64+54t_1^2+54t_2^2-24\sqrt{5}t_3+3\sqrt{10}t_1t_3+10t_3^2+48\sqrt{6}t_2-9\sqrt{30}t_2t_3$\\
		$\mathcal{B}_{4}$&$425+420t_1^2+420t_2^2-42\sqrt{5}t_3+5t_3^2+5\sqrt{2}t_1+3\sqrt{10}t_1t_3+345\sqrt{6}t_2-17\sqrt{30}t_2t_3$\\
		$\mathcal{B}_{3}$&$425+420t_1^2+420t_2^2+42\sqrt{5}t_3+5t_3^2+5\sqrt{2}t_1+3\sqrt{10}t_1t_3+345\sqrt{6}t_2+17\sqrt{30}t_2t_3$\\
		$\mathcal{B}_{2}$&$64+54t_1^2+54t_2^2+24\sqrt{5}t_3-3\sqrt{10}t_1t_3+10t_3^2+48\sqrt{6}t_2+9\sqrt{30}t_2t_3$\\
		$\mathcal{B}_{\overline{7}}$&$113+78t_1^2+78t_2^2+60\sqrt{5}t_3+35t_3^2+20\sqrt{2}t_1+12\sqrt{10}t_1t_3+76\sqrt{6}t_2+20\sqrt{30}t_2t_3$\\
		$\mathcal{B}_{\overline{1}}$&$-7-12t_1^2-12t_2^2+6\sqrt{5}t_3+5t_3^2-9\sqrt{6}t_2+\sqrt{30}t_2t_3+5\sqrt{2}t_1+3\sqrt{10}t_1t_3$\\
		\bottomrule
	\end{tabular}
	\normalsize
	\label{tab:isometric spheres}
\end{table}

\begin{table}[!htbp]
	\caption{The boundary arcs of the  octagon $\mathcal{E}$, which is the thick gray curve in  Figure \ref{figure:sidepair} of Section \ref{sec-3-orbifold}.}
	\centering
	\begin{tabular}{c|c|c}
		\toprule
		arc &Giraud disk& end points \\
		\midrule
		
		$\alpha_{1}$&$\partial_{\infty}(\mathcal{C}\cap \mathcal{B}_{\overline{0}})$  &$\partial_{\infty}(\mathcal{C}\cap \mathcal{B}_{\overline{0}}) \cap \partial_{\infty}\mathcal{B}_{\overline{1}}$, $\partial_{\infty}(\mathcal{C}\cap \mathcal{B}_{\overline{0}}) \cap \partial_{\infty}\mathcal{B}_{11}$ \\

		$\alpha_{2}$& $\partial_{\infty}(\mathcal{C}\cap \mathcal{B}_{11})$ &$\partial_{\infty}(\mathcal{C}\cap \mathcal{B}_{11}) \cap \partial_{\infty}\mathcal{B}_{\overline{0}}$, $\partial_{\infty}(\mathcal{C}\cap \mathcal{B}_{ 11}) \cap \partial_{\infty}\mathcal{B}_{5}$ \\

		$\alpha_{3}$&  $\partial_{\infty}(\mathcal{C}\cap \mathcal{B}_{5}))$ &$\partial_{\infty}(\mathcal{C}\cap \mathcal{B}_{5}) \cap \partial_{\infty}\mathcal{B}_{11}$, $\partial_{\infty}(\mathcal{C}\cap \mathcal{B}_{5}) \cap \partial_{\infty}\mathcal{B}_{4}$ \\

		$\alpha_{4}$& $\partial_{\infty}(\mathcal{C}\cap \mathcal{B}_{4})$  &$\partial_{\infty}(\mathcal{C}\cap \mathcal{B}_{4}) \cap \partial_{\infty}\mathcal{B}_{5}$, $\partial_{\infty}(\mathcal{C}\cap \mathcal{B}_{4}) \cap \partial_{\infty}\mathcal{B}_{3}$ \\

		$\alpha_{5}$&$\partial_{\infty}(\mathcal{C}\cap \mathcal{B}_{3})$ &$\partial_{\infty}(\mathcal{C}\cap \mathcal{B}_{3}) \cap \partial_{\infty}\mathcal{B}_{4}$, $\partial_{\infty}(\mathcal{C}\cap \mathcal{B}_{3}) \cap \partial_{\infty}\mathcal{B}_{2}$ \\
		
		$\alpha_{6}$&$\partial_{\infty}(\mathcal{C}\cap \mathcal{B}_{2})$ &$\partial_{\infty}(\mathcal{C}\cap \mathcal{B}_{2}) \cap \partial_{\infty}\mathcal{B}_{3}$, $\partial_{\infty}(\mathcal{C}\cap \mathcal{B}_{2}) \cap \partial_{\infty}\mathcal{B}_{\overline{7}}$ \\
		
		$\alpha_{7}$&$\partial_{\infty}(\mathcal{C}\cap \mathcal{B}_{\overline{7}})$  &$\partial_{\infty}(\mathcal{C}\cap \mathcal{B}_{\overline{7}}) \cap \partial_{\infty}\mathcal{B}_{2}$, $\partial_{\infty}(\mathcal{C}\cap \mathcal{B}_{\overline{7}}) \cap \partial_{\infty}\mathcal{B}_{\overline{1}}$  \\
		
		$\alpha_{8}$&$\partial_{\infty}(\mathcal{C}\cap \mathcal{B}_{\overline{1}})$  &$\partial_{\infty}(\mathcal{C}\cap \mathcal{B}_{\overline{1}}) \cap \partial_{\infty}\mathcal{B}_{\overline{7}}$, $\partial_{\infty}(\mathcal{C}\cap \mathcal{B}_{\overline{1}}) \cap \partial_{\infty}\mathcal{B}_{\overline{0}}$\\

		\bottomrule
	\end{tabular}
	\label{table:octagon}
\end{table}

The octagon $\mathcal{E}$  is bounded by  the eight segments on the intersections of $\partial_{\infty}\mathcal{C}$ with the closure of the  above eight bisectors. We denote by $\alpha_{1}$ $\alpha_{2}$, $\alpha_{3}$, $\alpha_{4}$, $\alpha_{5}$, $\alpha_{6}$, $\alpha_{7}$ and $\alpha_{8}$ the arcs on $\partial_{\infty}(\mathcal{C}\cap \mathcal{B}_{\overline{0}})$, $\partial_{\infty}(\mathcal{C}\cap \mathcal{B}_{11})$, $\partial_{\infty}(\mathcal{C}\cap \mathcal{B}_{5})$, $\partial_{\infty}(\mathcal{C}\cap \mathcal{B}_{4})$, $\partial_{\infty}(\mathcal{C}\cap \mathcal{B}_{3})$, $\partial_{\infty}(\mathcal{C}\cap \mathcal{B}_{2})$, $\partial_{\infty}(\mathcal{C}\cap \mathcal{B}_{\overline{7}})$ and  $\partial_{\infty}(\mathcal{C}\cap \mathcal{B}_{\overline{1}})$ respectively, see Table \ref{table:octagon} for the arcs and the boundaries of these eight arcs. From the equations in Table 	\ref{tab:isometric spheres} and the equation of $\partial_{\infty}\mathcal{C}$, one can deduce explicit
parametrizations for the segments of $\mathcal{E}$.

(1).   The resultant of  the equation of  $\partial_{\infty}(\mathcal{C}\cap \mathcal{B}_{\overline{0}})$ and $t_1^2+t_2^2+t_3^2-1$ with respect to $t_1$  has degree 2 in $t_2$. Using
the quadratic formula, we get
\begin{equation*}
t_2=\phi_1(t_3)=\frac{a_1(t_3)+\sqrt{b_1(t_3)}}{8(67+6\sqrt{5}t_3+15t_3^2)},
\end{equation*}
where
$$a_1(t)=-171\sqrt{6}-73\sqrt{30}t+123\sqrt{6}t^2+17\sqrt{30}t^3$$  and
$b_1(t)=8750-19500\sqrt{5}t+72250t^2-11400\sqrt{5}t^3-62750t^4+38580\sqrt{5}t^5-36810t^6.$

One then takes
\begin{equation*}
t_1=-\sqrt{1-\phi_1(t_3)^2-t_3^2}.
\end{equation*}

This parametrization is well defined for $t_3$ in the interval $[0,0.321084..]$ which corresponds to the segment on $\partial_{\infty}(\mathcal{C}\cap \mathcal{B}_{\overline{0}})$ of $\mathcal{E}$.
So we give a parametrization  for the arc  $\alpha_{1}$.

(2).  The segment $\alpha_{2}$ can be divided into two subsegments. We give a parametrization  for each one.
Let 
\begin{equation*}
t_2=\phi_2(t_3)=\frac{a_2(t_3)-\sqrt{b_2(t_3)}}{32(271-150\sqrt{5}t_3+105t_3^2)},
\end{equation*}
and
\begin{equation*}
t_2=\phi_2'(t_3)=\frac{a_2(t_3)+\sqrt{b_2(t_3)}}{32(271-150\sqrt{5}t_3+105t_3^2)},
\end{equation*}
where
$$a_2(t)=-955\sqrt{2}+873\sqrt{10}t-685\sqrt{2}t^2-129\sqrt{10}t^3,$$
and $b_2(t)=-2220150+9226020\sqrt{2}t-73068690t^2+60093240\sqrt{5}t^3-130836474t^4+27959460\sqrt{5}t^5-11466750t^6.$
Then we take
\begin{equation*}
\left[-\sqrt{1-\phi_2(t_3)^2-t_3^2}, \phi_2(t_3), t_3\right]
\end{equation*}  for $t_3$ in the interval $[0.242665..,0.270392..]$  for the first subsegment of $\alpha_{2}$  and 
\begin{equation*}
\left[-\sqrt{1-\phi_2'(t_3)^2-t_3^2}, \phi_2'(t_3), t_3\right]
\end{equation*}  for $t_3$ in the interval $[0.242665..,0.321084..]$  for the second subsegment of $\alpha_{2}$. Note that $\phi_2'(0.242665..)=\phi_2(0.242665..)$.

(3). For $\alpha_{3}$,  we get
\begin{equation*}
t_2=\phi_3(t_3)=\frac{a_3(t_3)-\sqrt{b_3(t_3)}}{12(192-72\sqrt{5}t_3+35t_3^2)},
\end{equation*}
where
\begin{eqnarray*}
	& a_3(t) =  -944\sqrt{6} + 369\sqrt{30}t + 172\sqrt{6} t^2-66\sqrt{30}t^3,  \\
	& b_3(t) = - 250t^2 + 1200\sqrt{5}t^3 - 1500\sqrt{5}t^4 +7680\sqrt{5}t^5-11140t^6.  \\
\end{eqnarray*}

Then we take
\begin{equation*}
t_1=-\sqrt{1-\phi_3(t_3)^2-t_3^2}.
\end{equation*}
This gives a parametrization  for $\alpha_{3}$ for $t_3$ in the interval $[0.171638..,0.270392..]$.

(4). For $\alpha_{4}$,  we get
\begin{equation*}
t_2=\phi_4(t_3)=\frac{a_4(t_3)-\sqrt{b_4(t_3)}}{8(17855-1758\sqrt{5}t_3+219t_3^2)},
\end{equation*}
where
\begin{eqnarray*}
	& a_4(t) =  -58305\sqrt{6} + 5771\sqrt{30}t + 27921\sqrt{6} t^2-1411\sqrt{30}t^3,  \\
	& b_4(t) = 350+1740\sqrt{5}t -17270t^2 + 42312\sqrt{5}t^3+17074t^4-306708\sqrt{5}t^5-651546t^6.  \\
\end{eqnarray*}

Then we take
\begin{equation*}
t_1=-\sqrt{1-\phi_4(t_3)^2-t_3^2}.
\end{equation*}
This gives a parametrization  for $\alpha_{4}$ for $t_3$ in the interval $[0,0.171638..]$.

(5). For $\alpha_{5}$,  we get
\begin{equation*}
t_2=\phi_5(t_3)=\frac{a_5(t_3)+\sqrt{b_5(t_3)}}{8(17855-1758\sqrt{5}t_3+219t_3^2)},
\end{equation*}
where
\begin{eqnarray*}
	& a_5(t) =  -58305\sqrt{6} - 5771\sqrt{30}t+ 27921\sqrt{6} t^2+1411\sqrt{30}t^3, \\
	& b_5(t) = 350-1740\sqrt{5}t -17270t^2 -42312\sqrt{5}t^3+17074t^4+306708\sqrt{5}t^5-651546t^6.  \\
\end{eqnarray*}

Then we take
\begin{equation*}
t_1=-\sqrt{1-\phi_5(t_3)^2-t_3^2}.
\end{equation*}
This gives a parametrization  for $\alpha_{5}$ for $t_3$ in the interval $[-0.171638..,0]$.

(6). For $\alpha_{6}$,  we get
\begin{equation*}
t_2=\phi_6(t_3)=\frac{a_6(t_3)+\sqrt{b_6(t_3)}}{12(192+72\sqrt{5}t_3+35t_3^2)},
\end{equation*}
where
\begin{eqnarray*}
	& a_6(t) =  -944\sqrt{6} - 369\sqrt{30}t + 172\sqrt{6} t^2+66\sqrt{30}t^3,  \\
	& b_6(t) = -250t_3^2 -1200\sqrt{5}t^3-9500t^4-7680\sqrt{5}t^5-11140t^6.  \\
\end{eqnarray*}

Then we take
\begin{equation*}
t_1=-\sqrt{1-\phi_6(t_3)^2-t_3^2}.
\end{equation*}
This gives a parametrization  for $\alpha_{6}$ for $t_3$ in the interval $[-0.270392..,-0.171638..]$.

(7). For $\alpha_{7}$,  there are two subsegments with one common end point. We give a parametrization  for each one.
Let
\begin{equation*}
t_2=\phi_7(t_3)=\frac{a_7(t_3)+\sqrt{b_7(t_3)}}{32(271+150\sqrt{5}t_3+105t_3^2)},
\end{equation*}
and 
\begin{equation*}
t_2=\phi_7'(t_3)=\frac{a_7(t_3)-\sqrt{b_7(t_3)}}{32(271+150\sqrt{5}t_3+105t_3^2)},
\end{equation*}
where
$$a_7(t)=-3629\sqrt{6}-2095\sqrt{30}t-683\sqrt{6}t^3+215\sqrt{30}t^3,$$
and $b_7(t)=-51250-247500\sqrt{5}t-2387750t^2-2452200\sqrt{5}t^3-7099550t^4-2180940\sqrt{5}t^5-1376010t^6.$

Then we take
\begin{equation*}
\left[-\sqrt{1-\phi_7(t_3)^2-t_3^2}, \phi_7(t_3),t_3\right]
\end{equation*}
for $t_3$ in the interval $[-0.321084..,-0.242665..]$ for the first part of $\alpha_{7}$  and 
\begin{equation*}
\left[-\sqrt{1-\phi_7'(t_3)^2-t_3^2}, \phi_7'(t_3),t_3\right]
\end{equation*}
for $t_3$ in the interval $[-0.270392, -0.242665]$ for the second part of $\alpha_{7}$.
Note that $\phi_7(-0.242665..)=\phi_7'(-0.242665..)$.

(8). For $\alpha_{8}$,  we get
\begin{equation*}
t_2=\phi_8(t_3)=\frac{a_8(t_3)-\sqrt{b_8(t_3)}}{8(67-6\sqrt{5}t_3+15t_3^2)},
\end{equation*}
where
$$a_8(t)=-171\sqrt{6}+73\sqrt{10}t+123\sqrt{6}t^2-17\sqrt{30}t^3,$$
and $b_8(t)=8750+19500\sqrt{5}t+72250t^2+11400\sqrt{5}t^3-62750t^4-38580\sqrt{5}t^5-36810t^6.$

Then we take
\begin{equation*}
t_1=-\sqrt{1-\phi_8(t_3)^2-t_3^2}.
\end{equation*}
This gives a parametrization  for $\alpha_{8}$ for $t_3$ in the interval $[-0.321084..,0]$.

\begin{figure}
	\centering
	\includegraphics[totalheight=3in]{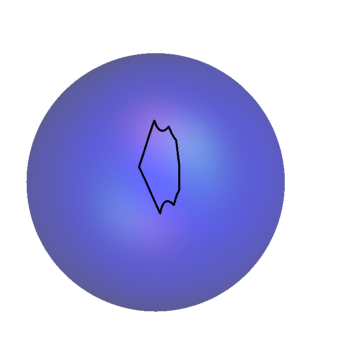}
	\caption{The cutting disk $\mathcal{E}$ on $\partial_{\infty}\mathcal{C}\simeq S^2$, which is an octagon.} \label{fig:cut-disk}
\end{figure}

The curve $\alpha=\bigcup_{i=1}^{8}\alpha_i$ on $\partial_{\infty}\mathcal{C}\simeq S^2$ bounds two disks, only one of which is completely contained in the half-sphere $t_1<0$. This is the cutting disk we need.
See Figure \ref{fig:cut-disk}.

We can show that the  boundary curve $\alpha$ of $\mathcal{E}$ is embedded in $\partial_{\infty} \mathcal{C}$ by solving  a system of equations. For example, the equations  $\partial_{\infty}(\mathcal{C}\cap \mathcal{B}_{\overline{0}})$,   $\partial_{\infty}(\mathcal{C}\cap \mathcal{B}_{11})$ and $t_1^2+t_2^2+t_3^2-1$ have two solutions
$$(-0.162508.., -0.933004.., 0.321084..),\quad( 0.0546295..,-0.856302.., 0.513578..).$$
 The point corresponding to the first solution is a vertex of $\mathcal{E}$.
From the parametrization for $\alpha_{i}$, we know that the point corresponding to the second solution is not on the curve $\alpha$.
Similar arguments apply to the intersections of the other pairs of arcs.

\begin{prop}\label{prop:octagon}  The topological octagon  $\mathcal{E}$ is properly contained in $T$.	
\end{prop}

\begin{proof}
	From the above construction, we see that points on the boundary of $\mathcal{E}$ are precisely on the bisectors we think they are on, see Figure \ref{fig:cut-disk}.
	
	We now want to check that the closure of the eight bisectors and $\partial_{\infty}\mathcal{C}$ have no unwanted
	extra intersection.
	
	For example, it is possible that $\mathcal{B}_{\overline{0}}$ may have a connected component contained in the interior of $\mathcal{E}$. In this case, the restriction to  $\partial_{\infty}\mathcal{C}$ of the equation for $\partial_{\infty}\mathcal{B}_{\overline{0}}\cap \partial_{\infty}\mathcal{C}$ would have a critical point in the interior of $\mathcal{E}$.
	
	By using  Lagrange multipliers, the critical points are the solutions of the system
	
	$$\left\{
	\begin{aligned}
	\nabla f&= \lambda \nabla g, \\
	g &=0, \\
	\end{aligned}
	\right.
	$$
	where $f$ is the equation for $\partial_{\infty}\mathcal{B}_{\overline{0}}\cap \partial_{\infty}\mathcal{C}$ and $g=t_1^2+t_2^2+t_3^2-1$.
	
	Then the system reads
	
	$$\left\{
	\begin{aligned}
	&-5\sqrt{2}+24t_1-2\lambda t_1+3\sqrt{10}t_3= 0, \\
	&9\sqrt{6}+24t_2-2\lambda t_2 +\sqrt{30}t_3=0, \\
	&6\sqrt{5}+3\sqrt{10}t_1+\sqrt{30}t_2-10t_3-2\lambda t_3 =0, \\
	&t_1^2+t_2^2+t_3^2-1 =0. \\
	\end{aligned}
	\right.
	$$

	Solving the system by standard Groebner basis techniques, we get two solutions which are given approximately
	
	$$(t_1,t_2,t_3)=(-0.173625, 0.942177, 0.286631),\,(0.308076, -0.341632, -0.887906).$$
	
	One can check that the points   correspond to these two solutions  are not in the interior of $\mathcal{E}$.  
	See Figure \ref{fig:critical-point-projection}, it is also clearly that the critical points of the equations are outside $\mathcal{E}$.

	\begin{figure}
		\centering
		\includegraphics[totalheight=2in]{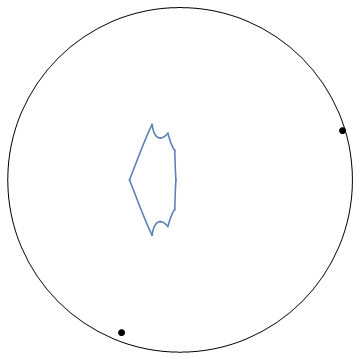}
		\caption{The projections of critical points and  the disk $\mathcal{E}$ onto the $(t_2,t_3)$-coordinates plane.} \label{fig:critical-point-projection}
	\end{figure}

	The analysis for the other intersections are similar, we omit the details.
This ends the proof of Proposition \ref{prop:octagon}.	
\end{proof}

\section{The 3-orbifold at infinity of $\rho(G_{6,3})$}\label{sec-3-orbifold}

 In this section,  based on  results in Section \ref{sec:combinatoricdirichlet},  we study the quotient of the domain of discontinuity under the action of the group  $\rho(G_{6,3})$ and identify the 3-orbifold  $\mathcal{O}$  at infinity of  ${\bf H}^2_{\mathbb C}/\rho(G_{6,3})$. We also  show $\mathcal{O}$ is a closed hyperbolic 3-orbifold. 

\subsection{The 3-orbifold $\mathcal{O}$ at infinity of $\rho(G_{6,3})$} \label{subsection:3orbifoldatinfty}

We have shown in Section  \ref{sec:combinatoricdirichlet} that $\partial_{\infty}D\cap \Omega$ is a solid torus $T$, and we have 
identified a simple closed curve in the boundary of $T$ which bounds a disk $\mathcal{E}$  in $T$.

\begin{figure}
	\begin{center}
		\begin{tikzpicture}
		\node at (0,0) {\includegraphics[width=11cm,height=9cm]{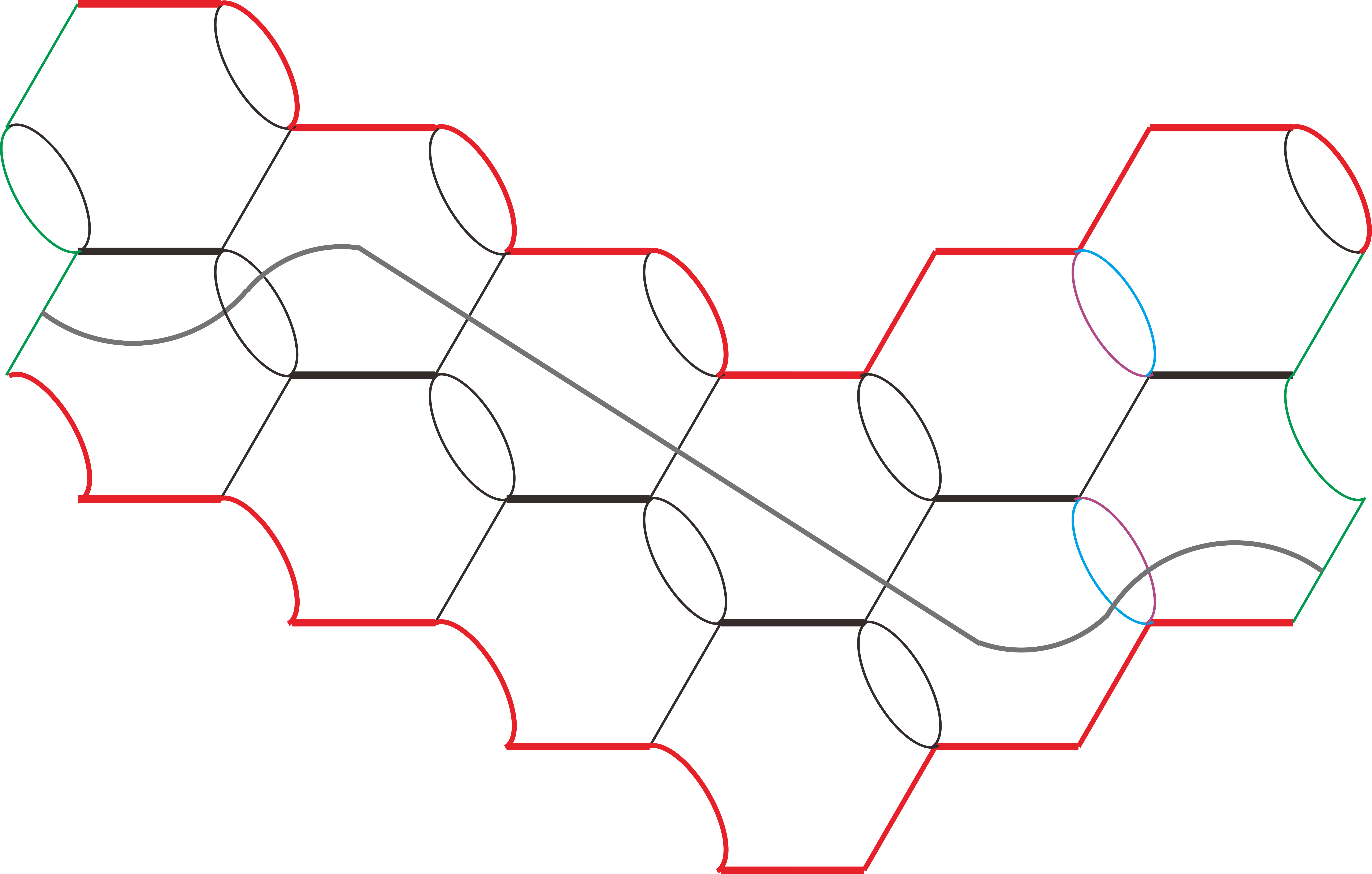}};
		
		\coordinate [label=left:\LARGE$0$] (S) at (-4.1,3.2);
		\coordinate [label=left:\LARGE$\bar{0}$] (S) at (-4.1,0.7);
		
			\coordinate [label=left:\LARGE$5$] (S) at (-2.3,1.9);
		\coordinate [label=left:\LARGE$\bar{5}$] (S) at (-2.3,-0.6);
		
			\coordinate [label=left:\LARGE$4$] (S) at (-0.5,0.6);
		\coordinate [label=left:\LARGE$\bar{4}$] (S) at (-0.5,-1.9);
		
			\coordinate [label=left:\LARGE$3$] (S) at (1.2,-0.7);
		\coordinate [label=left:\LARGE$\bar{3}$] (S) at (1.2,-3.2);

		\coordinate [label=left:\LARGE$2$] (S) at (2.9,-2);
		\coordinate [label=left:\LARGE$\bar{2}$] (S) at (2.9,0.7);

		\coordinate [label=left:\LARGE$1$] (S) at (4.6,2);
		\coordinate [label=left:\LARGE$\bar{1}$] (S) at (4.6,-0.7);

			\coordinate [label=left:\LARGE$11$] (S) at (-3.1,1.3);
		\coordinate [label=left:\LARGE$\bar{11}$] (S) at (-3.1,3.9);

		\coordinate [label=left:\LARGE$10$] (S) at (-1.4,0.1);
		\coordinate [label=left:\LARGE$\bar{10}$] (S) at (-1.4,2.6);

			\coordinate [label=left:\LARGE$9$] (S) at (0.2,-1.2);
		\coordinate [label=left:\LARGE$\bar{9}$] (S) at (0.2,1.3);
		
			\coordinate [label=left:\LARGE$8$] (S) at (1.95,-2.5);
		\coordinate [label=left:\LARGE$\bar{8}$] (S) at (1.95,0);
		
			\coordinate [label=left:\LARGE$7$] (S) at (3.6,1.3);
		\coordinate [label=left:\LARGE$\bar{7}$] (S) at (3.6,-1.2);
		
			\coordinate [label=left:\LARGE$6$] (S) at (-5.0,2.6);
		\coordinate [label=left:\LARGE$\bar{6}$] (S) at (5.35,2.6);

		\coordinate [label=left:\LARGE$A^{-1}$] (S) at (-2.5,-2.9);
		\coordinate [label=left:\LARGE$ A$] (S) at (0.7,2.9);

			\coordinate [label=left:$e_{1}$] (S) at (-4.0,2.1);
		\coordinate [label=left:$e_{2}$] (S) at (-2.25,0.8);
			\coordinate [label=left:$e_{3}$] (S) at (-0.5,-0.45);
				\coordinate [label=left:$e_{4}$] (S) at (1.3,-1.75);
					\coordinate [label=left:$e_{5}$] (S) at (3.0,-0.45);
						\coordinate [label=left:$e_{6}$] (S) at (4.7,0.8);
						
			\coordinate [label=left:$e_{7}$] (S) at (-4.0,-1.0);
				\coordinate [label=left:$e_{7}$] (S) at (-4.0,4.8);

		\coordinate [label=left:$e_{8}$] (S) at (-2.25,-2.2);
		\coordinate [label=left:$e_{8}$] (S) at (-2.25,3.4);

		\coordinate [label=left:$e_{9}$] (S) at (-0.5,-3.4);
			\coordinate [label=left:$e_{9}$] (S) at (-0.5,2.15);

		\coordinate [label=left:$e_{10}$] (S) at (1.3,-4.25);
		\coordinate [label=left:$e_{10}$] (S) at (1.3,0.9);
		
		\coordinate [label=left:$e_{11}$] (S) at (3.0,-3.4);
		\coordinate [label=left:$e_{11}$] (S) at (3.0,2.15);

		\coordinate [label=left:$e_{12}$] (S) at (4.7,-1.7);
		\coordinate [label=left:$e_{12}$] (S) at (4.7,3.4);

			\coordinate [label=left:$e_{13}$] (S) at (-5.1,3.95);
				\coordinate [label=left:$e_{13}$] (S) at (-5.3,2.5);
					\coordinate [label=left:$e_{13}$] (S) at (-4.8,-0.0);
					\coordinate [label=left:$e_{13}$] (S) at (5.25,1.5);
					\coordinate [label=left:$e_{13}$] (S) at (5.,0.1);
					\coordinate [label=left:$e_{13}$] (S) at (5.9,3.);

		\coordinate [label=left:$e_{14}$] (S) at (-4.4,2.9);
			\coordinate [label=left:$e_{14}$] (S) at (-5.12,1.3);
			\coordinate [label=left:$e_{14}$] (S) at (5.32,-0.99);
			\coordinate [label=left:$e_{14}$] (S) at (5.02,2.5);

			\coordinate [label=left:$e_{15}$] (S) at (-2.8,2.35);
		\coordinate [label=left:$e_{16}$] (S) at (-2.6,1.5);
		
			\coordinate [label=left:$e_{17}$] (S) at (-1.1,1.2);
		\coordinate [label=left:$e_{18}$] (S) at (-0.9,0.2);
		
			\coordinate [label=left:$e_{19}$] (S) at (0.7,0.1);
		\coordinate [label=left:$e_{20}$] (S) at (0.9,-1.2);
		
			\coordinate [label=left:$e_{21}$] (S) at (2.45,-1.2);
		\coordinate [label=left:$e_{22}$] (S) at (2.45,-2.2);

	\coordinate [label=left:$e_{23}$] (S) at (4.15,-2.5);
	\coordinate [label=left:$e_{24}$] (S) at (4.1,1.7);

		\end{tikzpicture}
	\end{center}
	\caption{The gluing pattern for   the 3-orbifold  $\mathcal{O}$. The edges $e_i$ and $e_{i+6}$ for $i=1, 2, 3, 4,5$ and $6$ are on the $\mathbb{C}$-circles fixed by the isometries $g_{0}, g_{5}, g_{4},  g_{3}, g_{2}$  and $g_{1}$ respectively.}
	\label{figure:sidepair}
\end{figure}

\begin{prop}\label{prop:idendify} The gray black curve and the red curve in  $\partial T$ in Figure \ref{figure:sidepair} are isotopic in the torus  $\partial T$.
\end{prop}

\begin{proof}This is trivial when we glue the sides in Figure \ref{figure:sidepair} to get a torus.
	\end{proof}

So now the red curve in Figure \ref{figure:sidepair} also bounds a disk in $T$.
Now we cut $T$ along this disk, we get a 3-ball $B$. Then the 3-orbifold $\mathcal{O}$ at infinity of $\rho(G_{6,3})$  is just the quotient space of $B$ with side pairings as in Figure \ref{figure:sidepair}. We now show this with more details.

Each hexagon labeled by $i$ (resp. $\bar{i}$) for $0 \leq i \leq 5$ in Figure  \ref{figure:sidepair} is part of the intersection of the bisector $\mathcal{B}_{i}$ (resp. $\mathcal{B}_{\bar{i}}$) and $\mathbb{S}^3$.  Each bigon labeled by $i$ (resp. $\bar{i}$) for $6 \leq i \leq 11$ in Figure  \ref{figure:sidepair} is part of the intersection of the bisector $\mathcal{B}_{i}$ (resp. $\mathcal{B}_{\bar{i}}$) and $\mathbb{S}^3$. 

Figure \ref{figure:sidepair} is a labeled disk.  Here we denote the union of three green   arcs (labeled by $e_{13}$, $e_{13}$ and $e_{14}$) in the hexagons and bigon labeled by $0$,  $\bar {0}$ and  $6$ by $\beta_{1}$, and we denote the union of three green   arcs (labeled by $e_{13}$, $e_{13}$ and $e_{14}$) in the hexagons  labeled by  $1$ and  $\bar{1}$ by $\beta_{2}$. Then, by gluing $\beta_{1}$ and $\beta_{2}$ together, we get an annulus in the boundary of the 3-ball $B$, such that the two thick red curves bound disjoint disks in the boundary of the 3-ball $B$. The upper thick red curve in Figure 	\ref{figure:sidepair} is denoted as $A$, and the lower thick red curve is denoted as $A^{-1}$. By abusing notation, we still use $A$ to denote the identity map gluing the red curve $A$ to the red curve $A^{-1}$, and $A^{-1}$ to denote the inverse of $A$. Then we glue the disks enclosed by $A$ and $A^{-1}$ getting the solid torus $T$. Moreover, if we glue together the green paths, and then the red paths in Figure 	\ref{figure:sidepair} (there is a twist when we glue the red circles, see the red edges labeled by $e_{13}$), we get a torus, which is the boundary of the solid torus $T$. 

For simplicity of notation, we write 
$$ g_i=\gamma_i \, \, {\rm and} \,\, g_{i+6}=\gamma_{i}\gamma_{i+1}^{-1},$$ where $i=0,1, 2, 3, 4, 5$.

Each $g^{-1}_{i}$ maps the hexagon (or bigon) labeled by $i$ to  the hexagon (or bigon) labeled by $\bar{i}$, and $A^{-1}$ maps the disk enclosed by $A$ to the disk enclosed by $A^{-1}$.  We now consider the actions of $\{g_{i}\}_{0 \leq i \leq 11}$ on the edges in Figure  \ref{figure:sidepair}:
\begin{itemize}
	\item
	
For $i=1$ (resp. $i=2$, $i=3$, $i=4$, $i=5$, $i=6$), the edge $e_{i}$ lies on the $\mathbb{C}$-circle fixed by $g_{0}$  (resp. $g_{5}$, $g_{4}$, $g_{3}$, $g_{2}$, $g_{1}$), so the element  $g_{0}$ (resp. $g_{5}$, $g_{4}$, $g_{3}$, $g_{2}$, $g_{1}$) fixes the edge $e_{1}$ (resp. $e_{2}$, $e_{3}$, $e_{4}$, $e_{5}$, $e_{6}$);

\item

 The $g^{-1}_{0}$-image of  the  arc labeled by $e_{14}$ in the boundary of the hexagon labeled by $0$ is the  green  arc labeled by $e_{14}$  in the boundary of the hexagon labeled by $\bar{0}$; 
 \item
 
 The $g^{-1}_{0}$-image of  the  green arc labeled by $e_{13}$ in the boundary of the hexagon labeled by $0$ is the red   arc labeled by $e_{13}$  in the boundary of the hexagon labeled by $\bar{0}$;
 
  \item
  From the above two actions, we can easily get the full action of $g^{-1}_{0}$ on the boundary arcs of the hexagon labeled by $0$;

 \item  Similarly, we have the action of $g_{i}$ for $i=1,2,3,4,5$ on these hexagons;
  
\item By taking sample points, we have  the $g^{-1}_{7}$-image  of  the blue arc in the boundary of the bigon labeled by $7$ (the intersection edge between the bigon labeled by $7$ and hexagon labeled by $1$) is the blue arc in the boundary of the bigon labeled by $\bar{7}$  (the intersection edge between the bigon labeled by $\bar{7}$ and hexagon labeled by $2$); The $g^{-1}_{7}$-image of the purple arc in the boundary of the bigon labeled by $7$ (the intersection edge between the bigon labeled by $7$ and hexagon labeled by $\bar{2}$)  is the purple arc in the boundary of the bigon labeled by $\bar{7}$ (the intersection edge between the bigon labeled by $\bar{7}$ and hexagon labeled by $\bar{1}$);

\item Similarly, we have the action of $g_{i}$ for $i=6,8,9,10,11$ on these bigons. 
\end{itemize}

From the side-pairings above,  we  get the  3-orbifold $\mathcal{O}$ at infinity of $\rho(G_{6,3})$.

\subsection{A presentation of  $\pi_{1}(\mathcal{O})$ } \label{subsection:fundamentalgroup}

From the  side-parings in Subsection \ref{subsection:3orbifoldatinfty},
 we will get a presentation of the fundamental group of  the 3-orbifold  $\mathcal{O}$ in this subsection.

 In Figure 	\ref{figure:sidepair}, we label the edge equivalent classes under the gluing pattern. But for the transparency of the figure, for some edge classes, we only label one of its representatives.  Then we get a presentation of $\pi_{1}(\mathcal{O})$ on thirteen generators $$g_0, g_1, g_2, \cdots, g_{11}, A$$ and twenty-four relations:
\begin{itemize}
	\item
	the edge $e_{1}$ (resp. $e_2$, $e_3$, $e_4$, $e_5$, $e_6$)  gives the relation $(g^{-1}_{0})^3=id$ (resp. $(g^{-1}_{5})^3=id$, $(g^{-1}_{4})^3=id$, $(g^{-1}_{3})^3=id$, $(g^{-1}_{2})^3=id$, $(g^{-1}_{1})^3=id$);
 	\item
 	
 	the edge $e_{7}$ (resp. $e_8$, $e_9$, $e_{10}$, $e_{12}$)  gives the relation $(Ag^{-1}_{0})^3=id$ (resp. $(Ag^{-1}_{5})^3=id$, $(Ag^{-1}_{4})^3=id$, $(Ag^{-1}_{3})^3=id$, $(Ag^{-1}_{1})^3=id$); and the edge $e_{11}$ gives the relation $(Ag_{2})^3=id$;
 	 	
 %	 the edge $e_{i+6}$ give the relation $(Ag^{-1}_{i-1})^3=id$ for $i=1,2,4,5,6$, $e_{11}$ %give the relation $(Ag_{2})^3=id$;
 	 	\item  there are twelve more relations given in Table \ref{table:relation}.
\end {itemize}

\begin{table}[!htbp]
	\caption{Some cycle relations  of the 3-orbifold  $\mathcal{O}$.}
	\centering
	\begin{tabular}{c|c}
		\toprule
		\textbf{Ridge} & \textbf{Cycle relation  } \\
		\midrule
		$e_{13}$  & $g_1g_6Ag^{-1}_{0}$ \\

		$e_{14}$  & $g_0g^{-1}_{1}g^{-1}_6$ \\
		
			$e_{15}$  &  $g_0g_{11}Ag^{-1}_{5}$ \\

		$e_{16}$  & $g_5g^{-1}_{0}g^{-1}_{11}$ \\
		
			$e_{17}$  & $g_5g_{10}Ag^{-1}_{4}$ \\

		$e_{18}$  & $g_4g^{-1}_{5}g^{-1}_{10}$ \\
		
			$e_{19}$  & $g_4g_9Ag^{-1}_{3}$ \\
		
		$e_{20}$  & $g_3g^{-1}_{4}g^{-1}_9$ \\

		$e_{21}$  & $g_{3}g_8g^{-1}_2$ \\
			$e_{22}$  & $g_8g_3A^{-1}g^{-1}_{2}$ \\
		
			$e_{23}$  & $g_2g_7g^{-1}_{1}A$ \\

		$e_{24}$  &$g_7g_{2}g^{-1}_1$ \\

		\bottomrule
	\end{tabular}
	\label{table:relation}
\end{table}

Then we can cancel the generators $g_6$, $g_7$, $\cdots$, $g_{11}$ to give a presentation of $\pi_{1}(\mathcal{O})$:
$$\pi_{1}(\mathcal{O})=  \left\langle g_0, \cdots,  g_5,A \Bigg| \begin{matrix}  g_{i}^{3}=id, \quad 
(Ag^{-1}_{0})^3=(Ag^{-1}_{1})^3=(Ag^{-1}_{3})^3=id, \\ (Ag^{-1}_{4})^3= (Ag^{-1}_{5})^3= (Ag_{2})^3=id, \\
Ag^{-1}_0g_1g_0g^{-1}_1=Ag^{-1}_5g_0g_5g^{-1}_0= Ag^{-1}_4g_5g_4g^{-1}_5=\\ Ag^{-1}_3g_4g_3g^{-1}_4
= Ag^{-1}_3g^{-1}_2g_3g_2=Ag_2g_1g^{-1}_2g^{-1}_1=id\end{matrix}\right\rangle 
.$$

We should remark here the presentation of $\pi_{1}(\mathcal{O})$   is not symmetric with the generators above, and moreover if  we add the relation $A=id$ to $\pi_{1}(\mathcal{O})$, then we get the original group $G_{6,3}$.

Using the  Magma Calculator available at \cite{Mag}, we can simplify the presentation to get 
\begin{equation}\label{eq:groupviaDirichlet}
\pi_{1}(\mathcal{O})=  \left\langle s_1, \cdots,  s_6\Bigg| \begin{matrix}  s_{i}^{3}=id, \quad  s_2 s^{-1}_1 s^{-1}_2  s_1  s^{-1}_6 s_1  s_6  s^{-1}_1= id, \\
s^{-1} s_2  s_1 s_3  s^{-1}_2  s^{-1}_3=id, \quad s_5  s_4  s^{-1}_5 s^{-1}_3 s^{-1}_4 s_3=id,\\

s_6 s^{-1}_5 s^{-1}_6 s_5 s^{-1}_1 s_2 s_1  s^{-1}_2 =id,\quad  s_2  s^{-1}_1 s^{-1}_2 s_1  s^{-1}_4 s_5 s_4 s^{-1}_5 = id
\end{matrix}\right\rangle.
\end{equation}

Changing $u_{3}=s_2s_3s^{-1}_2$,  $s_3=s^{-1}_{2}u_{3}s_2$, $u_{i}=s_{i}$ for $i=1,2,4,5,6$,  we can rewrite the presentation  as 
\begin{equation}\label{eq:groupviaDirichletlast}
	\pi_{1}(\mathcal{O})=  \left\langle u_1, u_2,u_3,u_4,u_5,  u_6\Bigg| \begin{matrix}  u_{i}^{3}=id,  \quad 
		u_3u^{-1}_2u^{-1}_3u_2=	u_2u^{-1}_1u^{-1}_2u_1=id,\\
			u_1u^{-1}_6u^{-1}_1u_6=	u_6u^{-1}_5u^{-1}_6u_5=id,\\
			u_5u^{-1}_4u^{-1}_5u_4= u_5u_4u^{-1}_5u^{-1}_2u^{-1}_3u_2u^{-1}_4u^{-1}_2u_3u_2=id
		 \end{matrix}\right\rangle.
\end{equation}

\subsection{A chain link orbifold}
\begin{figure}
	\begin{center}
		\begin{tikzpicture}
		\node at (0,0) {\includegraphics[width=9cm,height=9cm]{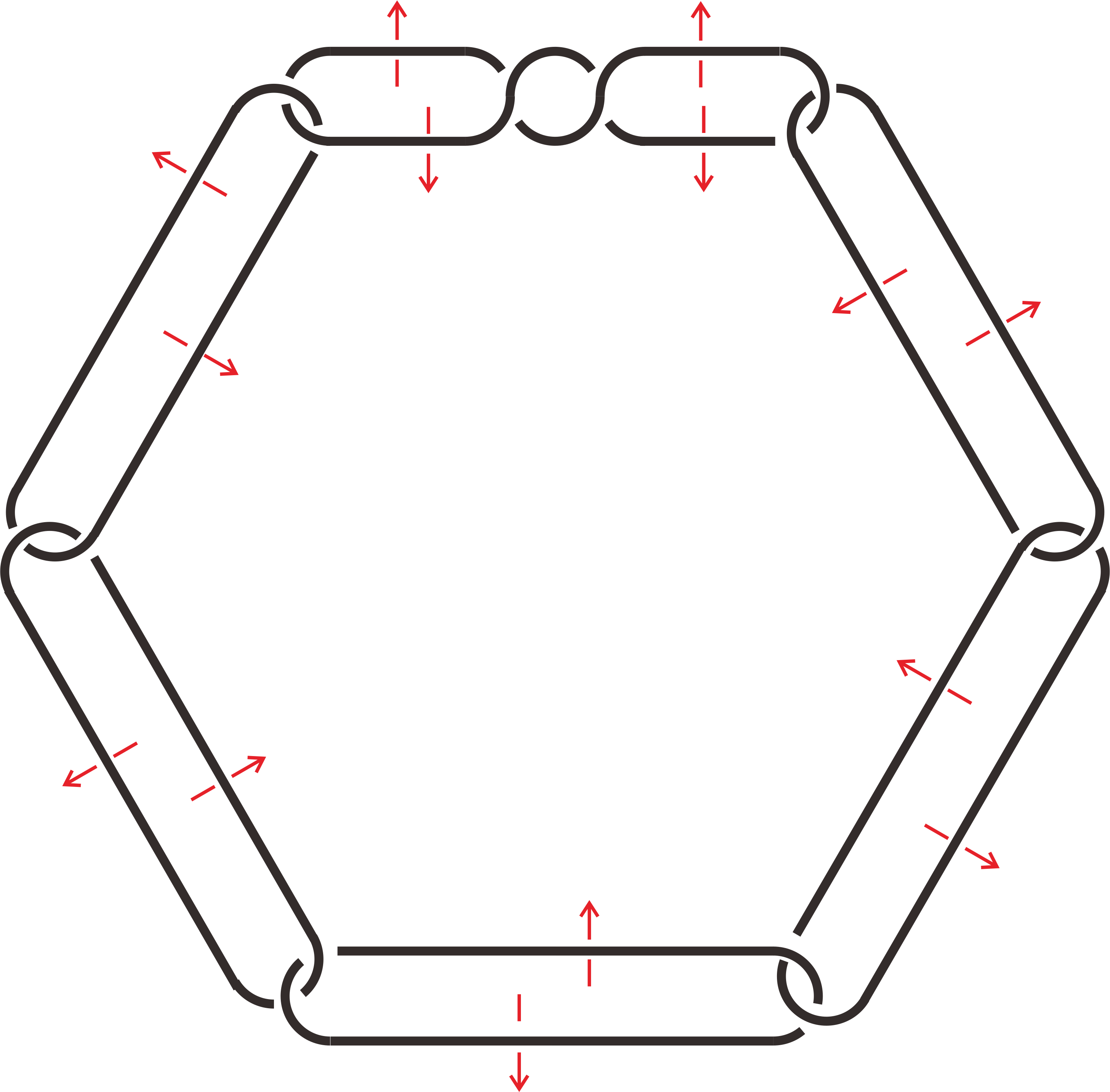}};

		\coordinate [label=left:$y_{6}$] (S) at (-1.3,4.3);
		\coordinate [label=left:$z_{6}$] (S) at (-1.05,3);

		\coordinate [label=left:$y_{5}$] (S) at (1.2,4.3);
		\coordinate [label=left:$z_{5}$] (S) at (1.2,3);

			\coordinate [label=left:$y_{4}$] (S) at (4.2,2.3);
		\coordinate [label=left:$z_{4}$] (S) at (2.5,1.7);

	\coordinate [label=left:$y_{3}$] (S) at (4.2,-2.3);
	\coordinate [label=left:$z_{3}$] (S) at (2.8,-1.2);
	
		\coordinate [label=left:$y_{2}$] (S) at (-0.4,-4.3);
	\coordinate [label=left:$z_{2}$] (S) at (0.2,-3);
	
			\coordinate [label=left:$y_{1}$] (S) at (-3.4,-2.3);
		\coordinate [label=left:$z_{1}$] (S) at (-1.9,-2.0);
		
			\coordinate [label=left:$y_{0}$] (S) at (-3.2,2.9);
		\coordinate [label=left:$z_{0}$] (S) at (-2.0,1.7);

		\end{tikzpicture}
	\end{center}
	\caption{The Chain link $C(6,-2)$.}
	\label{figure:chainlink}
\end{figure}

Consider the link in Figure 	\ref{figure:chainlink}, it is called the chain link $C(6,-2)$ in \cite{NReid} (in fact, our link here is the mirror image of the  link $C(6,-2)$ in \cite{NReid}, but this does not matter, since they have  homeomorphic complements.)

Then from the standard Wirtinger presentation of the fundamental group of a link in the 3-sphere \cite{Rolfsen}, we see that $ \pi_{1}(\mathbb{S}^3-C(6,-2))$ is a group on fourteen generators $$  y_0,\ \cdots,  \ y_6, \ z_0, \ \cdots,  \ z_6$$ and fourteen  relations in Table 	\ref{table:relationofchain}.
\begin{table}[!htbp]
	\caption{ Relations  of $ \pi_{1}(\mathbb{S}^3-C(6,-2))$.}
	\centering
	\begin{tabular}{c|c}
		\toprule
				crossing &relation \\
				\midrule
				
					1&	$z_6y_0y^{-1}_{6}y^{-1}_{0}$ \\
				2&	$z_6z_{0}z^{-1}_6y^{-1}_{0}$ \\

				3&  $z_0y_1y^{-1}_{0}y^{-1}_{1}$ \\
				4& $z_0z_{1}z^{-1}_0y^{-1}_{1}$ \\

	5&	$z_1y_2y^{-1}_{1}y^{-1}_{2}$ \\
		6&	$z_1z_{2}z^{-1}_1y^{-1}_{2}$ \\
		
			7&	$z_2y_3y^{-1}_{2}y^{-1}_{3}$ \\
		8&	$z_2z_{3}z^{-1}_2y^{-1}_{3}$ \\
		
		9&	$z_3y_4y^{-1}_{3}y^{-1}_{4}$ \\
		10&	$z_3z_{4}z^{-1}_3y^{-1}_{4}$ \\
		
			11&	$z_4y_5y^{-1}_{4}y^{-1}_{5}$ \\
		12&	$z_4z_{5}z^{-1}_4y^{-1}_{5}$ \\
		
			13&	$z_5z_6z^{-1}_{6}y^{-1}_{6}$ \\
		14&	$y_5z_{6}y^{-1}_5z^{-1}_{5}$ \\

		\bottomrule
	\end{tabular}
	\label{table:relationofchain}
\end{table}

Here the relations $z_5z_6z^{-1}_{6}y^{-1}_{6}$ and 
$y_5z_{6}y^{-1}_5z^{-1}_{5}$ correspond to the two crossings in the twist region  of Figure 	\ref{figure:chainlink}, so they have different forms from those of the  others.

We consider the 3-orbifold  $\mathcal{L}$ with underlying space the 3-sphere  and whose singular locus is the   $\mathbb{Z}_3$-coned chain-link $C(6,-2)$. Then its fundamental group is just add the relations $y_{i}^3=z_{i}^3=id$ for $i=0,1,2,3,4,5,6$ to the presentation $\pi_{1}(\mathbb{S}^3-C(6,-2))$ above.

Using Magma \cite{Mag}, we can simplify the  presentation of  $\pi_{1}(\mathcal{L})$ to get
$$\pi_{1}(\mathcal{L})=  \left\langle t_1, \cdots,  t_6\Bigg| \begin{matrix}  t_{i}^{3}=id,  \quad 
t^{-1}_5t_6t_5t^{-1}_6t_4t^{-1}_3t^{-1}_4t_3=t^{-1}_5t_6t_5t^{-1}_6t_2t^{-1}_1t^{-1}_2t_1=id,
\\
t^{-1}_2t_3t_2t^{-1}_3t_6t^{-1}_5t^{-1}_6t_5=t^{-1}_5t_6t_5t^{-1}_6t_5t^{-1}_4t^{-1}_5t_4=id,\\ 
t^{-1}_1t_2t_1t^{-1}_2t_3t^{-1}_2t^{-1}_3t_2=t_2t_1t^{-1}_2t^{-1}_5t^{-1}_6t_5t^{-1}_1t^{-1}_5t_6t_5=id\\
\end{matrix}\right\rangle 
.$$

We can rewrite it as 
\begin{equation}\label{eq:groupvialink}
\pi_{1}(\mathcal{L})=  \left\langle t_1, \cdots,  t_6\Bigg| \begin{matrix}  t_{i}^{3}=id,  \quad 
t_6t^{-1}_5t^{-1}_6t_5=t_5t^{-1}_4t^{-1}_5t_4=id,\\
t_4t^{-1}_3t^{-1}_4t_3=t_3t^{-1}_2t^{-1}_3t_2=id,\\
t_2t^{-1}_1t^{-1}_2t_1=
t_2t_1t^{-1}_2t^{-1}_5t^{-1}_6t_5t^{-1}_1t^{-1}_5t_6t_5=id
\end{matrix}\right\rangle 
.
\end{equation}

Now the map $f:\pi_{1}(\mathcal{L}) \rightarrow\pi_{1}(\mathcal{O})$, 
\begin{equation}\label{eq:iso}
 \begin{matrix} 
t_1 \rightarrow u_4,\\
t_2 \rightarrow u_5,\\
t_3\rightarrow u_6,\\
t_4 \rightarrow u_1,\\
t_5 \rightarrow u_2,\\
t_6 \rightarrow u_3,\\
\end{matrix}
\end{equation}
is  an isomorphism between $\pi_{1}(\mathcal{L})$ in the presentation (\ref{eq:groupvialink}) and $\pi_{1}(\mathcal{O})$ in the presentaition (\ref{eq:groupviaDirichletlast}).

Now we take a torsion free finite index subgroup of $\pi_{1}(\mathcal{L})=\pi_{1}(\mathcal{O})$, which corresponds  to finite coves $\widetilde{\mathcal{L}}\rightarrow \mathcal{L}$  and $\widetilde{\mathcal{O}}\rightarrow \mathcal{O}$. We will show in Subsection \ref{subsection:hyper}  that $\mathcal{L}$ is a  closed  hyperbolic 3-orbifold. So $\widetilde{\mathcal{L}}$ is a closed hyperbolic  3-manifold.
 Then, by the prime decompositions of 3-manifolds \cite{Hempel}, $\widetilde{\mathcal{O}}$ is the connected sum of $\widetilde{\mathcal{L}}$ with  $\mathcal{N}$, where
 $\mathcal{N}$  is a closed 3-manifold with trivial fundamental group. By the solution of the
Poincar\'e Conjecture, then   $\mathcal{N}$  is the 3-sphere, so $\widetilde{\mathcal{L}}$ is homeomorphic to  $\widetilde{\mathcal{O}}$. This in turn implies $\mathcal{L}$ is homeomorphic to  $\mathcal{O}$.
This finishes the proof of the first part of Theorem \ref{thm:main}.
\subsection{Hyperbolicity of  the 3-orbifold $\mathcal{O}$} \label{subsection:hyper}

	 We show the orbifold $\mathcal{O}$ is hyperbolic by the arguments from Cooper-Hodgson-Kerckhoff \cite{CHK} and    Dunbar \cite{dunbar}. From Theoren 1.25 of \cite{CHK}, we see $\mathcal{O}$  is a geometric 3-orbifold. That is, the 3-orbifold with underlying space $\mathbb{S}^3$ and $\mathbb{Z}_3$-coned singular set  on  $C(6,-2)$ is a geometric 3-orbifold. In page 81 of \cite{dunbar}, there is classification of non-hyperbolic 3-orbifold with underlying space $\mathbb{S}^3$.   There is no  orbifold with singular set a $\mathbb{Z}_3$-coned six components link    in the list.  $\mathcal{O}$ is hyperbolic, since it is not one of the exceptions.   We thank one of the referees for pointing this out to us.

	This completes the proof of Theorem \ref{thm:main}.

%\begin{figure}
%	\begin{center}
%		\begin{tikzpicture}
%		\node at (0,0) {\includegraphics[width=6cm,height=7cm]{surgerywhitehead}};
%		
%			\coordinate [label=left:$\frac{6}{-1}$] (S) at %(-3.3,0.3);
%			\coordinate [label=left:$\frac{3}{0}$] (S) at (3.25,2.7);
%	\coordinate [label=left:$y_{6}$] (S) at (-3.3,0.3);
%	\coordinate [label=left:$z_{6}$] (S) at (3.25,2.7);		
%\end{tikzpicture}
%	\end{center}
%\caption{By   performing $\frac{6}{-1}$  Dehn filling on one component of the Whitehead link, and  performing $\frac{3}{0}$ Dehn filling on the other  component, we get  the 3-orbifold $\mathcal{Q}$.}
%\label{figure:whiteheadlink}
%\end{figure}

%By   performing$\frac{6}{-1}$  $(6,-1)$ Dehn filling on one component of the Whitehead link, and  performing $\frac{6}{-1}$ $(3,0)$ Dehn filling on the other  component, we get performing the 3-orbifold $\mathcal{Q}$.

\end{document}